\documentclass[aap,MSNbibl,mtplusscr,seceqn,citesort,dvips]{arximspdf}
\usepackage{mathrsfs}
\usepackage[vtexbibl]{}
\usepackage{graphicx}

\doi{10.1214/09-AAP665}
\volume{20}
\issue{4}
\pubyear{2010}
\firstpage{1219}
\lastpage{1252}

\makeatletter

\newtheorem{theorem}{Theorem}[section]
\newtheorem{lemma}{Lemma}[section]
\newtheorem{proposition}{Proposition}[section]

\newproclaim{definition}{Definition}[section]
\newproclaim{example}{Example}[section]
\newproclaim{remark}{Remark}[section]

\def\eqref#1{(\ref{#1})}
\makeatother

\begin{document}
\begin{frontmatter}

\title{Convergence of complex multiplicative cascades}
\runtitle{Convergence of complex multiplicative cascades}

\begin{aug}
\author[A]{\fnms{Julien} \snm{Barral}\ead[label=e1]{barral@math.univ-paris13.fr}\corref{}},
\author[B]{\fnms{Xiong} \snm{Jin}\ead[label=e2]{Xiong.Jin@inria.fr}}
\and
\author[C]{\fnms{Beno\^\i t} \snm{Mandelbrot}\ead[label=e3]{Benoit.Mandelbrot@yale.edu}}
\affiliation{Universit\'e Paris 13, INRIA and Yale University}
\address[A]{J. Barral\\ LAGA (UMR 7539), Institut Galil\'ee \\
Universit\'{e} Paris 13\\ 99 avenue Jean-Baptiste Cl\'ement \\ 93430
Villetaneuse\\ France \\\printead{e1}}
\address[B]{X. Jin\\ INRIA\\ Domaine de Voluceau \\78153 Le Chesnay cedex \\
France\\ \printead{e2}}
\address[C]{B. Mandelbrot\\Department of Mathematics \\Yale University \\
New Haven, Connecticut 06520 \\USA\\ \printead{e3}}
\runauthor{J. Barral, X. Jin and B. Mandelbrot}
\end{aug}

\received{\smonth{2} \syear{2009}}
\revised{\smonth{11} \syear{2009}}


%
\begin{abstract}
The familiar cascade measures are sequences of random positive
measures obtained on $[0,1]$ via $b$-adic independent cascades. To
generalize them, this paper allows the random weights invoked in the
cascades to take real or complex values. This yields sequences of
random functions whose possible strong or weak limits are natural
candidates for modeling multifractal phenomena. Their asymptotic
behavior is investigated, yielding a sufficient condition for almost
sure uniform convergence to nontrivial statistically self-similar
limits. Is the limit function a monofractal function in multifractal
time? General sufficient conditions are given under which such is the
case, as well as examples for which no natural time change can be used.
In most cases when the sufficient
condition for convergence does not hold, we show that either the
limit is 0 or the sequence diverges almost surely. In the later case,
a functional central limit theorem holds, under some conditions. It
provides a natural normalization making the sequence converge in law
to a standard Brownian motion in multifractal time.

\end{abstract}



\end{frontmatter}
%

\section{Introduction}\label{intro}

\subsection{Foreword about H\"{o}lder singularity and multifractal functions}\label{1.1}
Multifractal analysis is a natural framework to
describe statistically and geometrically the heterogeneity in the
distribution at small scales of the H\"older singularities of a given
locally bounded function or signal $F\dvtx I\mapsto\mathbb{R}$ (or $\mathbb
{C}$) where $I$ is a bounded interval. One possible way to define the
H\"older singularity at a given point $t$ is by measuring the
asymptotic behavior of the oscillation of $f$ around $t$ thanks to the exponent
\[
h_F(t)=\liminf_{r\to0^+} \frac{\log\operatorname{Osc}_F([t-r,t+r])}{\log(r)}
\]
 or
 \[
 h_F(t)=\liminf_{n\to\infty}
\frac{\log_2\operatorname{Osc}_F(I_n(t))}{-n},
\]
where $I_n(t)$ is the dyadic interval of length $2^{-n}$ containing $t$
and $\operatorname{Osc}_F(J)=\sup_{s,t\in J}|F(t)-F(s)|$. Then
the multifractal analysis of $F$ consists in classifying the points
with respect to their H\"older singularity.
The singularity spectrum of $F$ is the Hausdorff dimension of the H\"
older singularities level sets
\[
E_F(h)=\{t\in I\dvtx  h_F(t)=h\},\qquad h\ge0.
\]
One says that $F$ is monofractal if there exists a unique $h\ge0$ such
that $E_F(h)\neq\varnothing$. Otherwise $F$ is multifractal.
Alternatively, one can also compute free energy
functions like
%
\begin{equation}\label{tauF}
\tau_F(q)=\liminf_{n\to\infty}\frac{-1}{n}\log_2\sum_{I\in\mathcal
{G}_n,\operatorname{Osc}_F(I)\neq0}\operatorname{Osc}_F(I)^q
\end{equation}
and say that $F$ is monofractal if $\tau_F$ is linear. Let us
mention that one always has $\dim E_F(h)\le
\tau_F^*(h)=\inf_{q\in\mathbb{R}}hq-\tau_F(q)$, and one says that the
multifractal formalism holds at $h$ if these inequalities are
equalities (see \mbox{\cite{JAFFJMP,Riedi}} for
instance). Also, one says that $F$ is monofractal in the strong sense
if $ \lim_{r\to0^+} \frac{\log{\rm
Osc}_{F}([t-r,t+r])}{\log(r)}$ exists everywhere and is independent of
$t$. Examples of functions that satisfy this property are the
Weierstrass functions,
Brownian motion and self-affine functions \cite{Urb1,Heurteaux}.

When $F$ is a continuous multifractal function possessing some scaling
invariance property, it may happen that it possesses the remarkable
property to be decomposable as a (often strongly) monofractal function
$B$ and an increasing multifractal time change $G$ such that $F=B\circ
G$. In the known cases, there exists $\beta>1$ such that $\tau_F(\beta
)=0$, and the function $G$ is such that, roughly speaking,
$|G(I)|\approx\operatorname{Osc}_F(I)^{\beta}$ for $I$ small enough in $\bigcup
_{n\ge1}\mathcal{G}_n$, so that $\tau_{G}(q)=\tau_F(\beta q)$ (\cite
{Mandfin,Sadv}).

\subsection{Some methods that build multifractal processes} Our goal is
to build new multifractal stochastic processes. Since \cite{M1,M2,M3},
one of the main
approaches is to construct singular statistically
self-similar measures $\mu$ on $[0,1]$ (or
$\mathbb{R}_+$). These measures are obtained as almost sure weak limits
of absolutely
continuous measure-valued martingales $(\mu_n)_{n\ge1}$ whose
densities $(Q_n)_{n\ge1}$ are martingales generated by
multiplicative processes (see \cite{M2,M3,BM2,Mann,BacryMuzy,BMS2}).
These objects have been used to construct nonmonotonic multifractal
stochastic processes as follows: (a) by starting with fractional
Brownian motions or stable L\'evy processes
$X$, and an independent measure $\mu$, and performing the
multifractal time change $X(\mu([0,t]))$
\cite{Mandfin,BacryMuzy,Riedi,CAR,BS}; (b) by integrating the
density $Q_n$ with respect to the Brownian motion and letting $n$
tend to $\infty$ \cite{BacryMuzy,CAR}; (c) by using $\mu$ to specify
the covariance of
some Gaussian processes \cite{Ludena}; (d) by considering random wavelet
series whose coefficients are built from the multifractal measure
$\mu$ \cite{ABM,BSw}. Such processes are geared to modeling signals
possessing a wildly varying
local H\"older regularity as well as scaling invariance properties.
Many such signals come from physical or social intermittent
phenomena like turbulence \cite{M3,FrPa}, spatial rainfall
\cite{Gupta}, human heart rate \cite{Stanley,Meyer}, Internet
traffic \cite{RRJLV,Willinger} and stock exchanges prices
\cite{Mandfin,BacryMuzy2}.

As background, standard statistically self-similar measures consist in the
1-dimensional $b$-adic canonical cascades constructed on the interval
$[0,1]$ as follows.
Fix an integer $b\ge2$. The $b$-adic closed subintervals of $[0,1]$
are naturally encoded by the nodes of the tree $\mathscr{A}^*=\bigcup
_{n\ge0}\{0,\dots,b-1\}^n$ with the convention
that $\{0,\dots,b-1\}^0$ contains the root of $\mathscr{A}^*$ denoted
$\varnothing$. To each element $w$ of $\mathscr{A}^*$,
we associate a nonnegative random weight $W(w)$, these
weights being independent and identically distributed with a random
variable $W$ such that
$\mathbb{E}(W)=1$. A~sequence of random densities $(Q_{n})_{n\geq1}$
is then obtained as follows: Let $I$ be the semi-open to the right
$b$-adic interval encoded by
the node $w=w_1w_2\cdots w_n$, that is,
$I=[\sum_{k=1}^nw_kb^{-k},b^{-n}+\sum_{k=1}^nw_kb^{-k}]$. Then
\[
Q_n(t)=W(w_1)W(w_1w_2)\cdots
W(w_1w_2\cdots w_n)\qquad \mbox{for $t\in I$}.
\]
Consider the sequence of measures $(\mu_n)_{n\ge1}$ whose
densities with respect to the Lebesgue measure are given by
$(Q_{n})_{n\geq1}$. This is the measure-valued
martingale (with respect to the natural filtration it generates)
introduced in \cite{M2,M3}.

The study of these martingales and their limits has led to numerous
mathematical developments
in the theories of probability and geometric measure \cite
{KP,DL,KAIHP,K2,K3,K4,Gu,CK,HoWa,Mol,Fal,AP,B1,B2,LR,OW,Fan5,Liu2002,LRR,BS1,WaWi1,WaWi}.
These objects, as well as the other statistically self-similar measures
mentioned above, are special
examples of a general model of positive measure-valued martingale,
namely the ``$T$-martingales'' developed in
\cite{K2,K3}, which make rigorous the construction and results of the
seminal work \cite{M1} on log-normal multiplicative chaos.

A natural alternative to the preceding constructions consists of
allowing the
$[0,1]$-martingales $(Q_n)_{n\ge1}$ to take real or complex values and
consider the continuous functions-valued
martingale
\[
\biggl(F_n\dvtx t\in[0,1]\mapsto\int_0^tQ_n(s)\, \mathrm{d}s \biggr)_{n\ge1}.
\]
The companion paper \cite{BJMpartI}
exhibits a class of such $[0,1]$-martingales for which we can prove a
general uniform convergence theorem
to a nontrivial random function. As a consequence, we construct
natural extensions of random functions of
the now familiar statistically self-similar measures, namely canonical
$b$-adic cascades \cite{M2,M3}, compound
Poisson cascades \cite{BM2,B4} and infinitely divisible cascades \cite
{BacryMuzy,CAR}.

\subsection{Further results for complex $b$-adic cascades} This paper
deepens the study of the convergence
properties of the complex extension of 1-dimensional canonical $b$-adic
cascades. In particular, we improve
the convergence result obtained in \cite{BJMpartI}, and proceed beyond
the results obtained in \cite{BM} in
the special case where $W$ is real-valued and constant in absolute
value. In this case, two possibilities exist as $n$ tends to $\infty$.
$(F_n)_{n\ge1}$ may
converge almost surely uniformly to a monofractal process sharing some
fractal properties
with a fractional Brownian motion of exponent $H\in(1/2,1)$. If not,
$F_n$ is not bounded and the following new functional
central limit theorem holds: as $n$ tends to $\infty$, $F_n/\sqrt
{\mathbb{E}(F_n(1)^2)}$ converges in law to the
restriction to $[0,1]$ of the standard Brownian motion. This last
result raises a question. Does
$F_n/\sqrt{\mathbb{E}(F_n(1)^2)}$ have a weak limit in the general case
when $F_n(1)$ is not bounded in $L^2$ norm?
In case of weak convergence, it remains to describe the nature of the
limit multifractal process and compare it with other
models of random multifractal functions.

On the asymptotic behavior of $(F_n)_{n\ge1}$, our main results are of
the following nature. We obtain a
condition on the moments of $W$ that suffices for the uniform
convergence of $F_n$---almost surely and in $L^p$
norm ($p>1$)---to a nontrivial limit (Theorem \ref{conv-1}). When
this sufficient condition does not hold,
we show that---in most of the cases---either $F_n$ converges
uniformly to 0, or $F_n$ diverges almost surely in
$\mathcal{C}([0,1])$, namely in the space of complex-valued continuous
functions over $[0,1]$ (Theorems \ref{div} and \ref{div-caspart}); it
is worth mentioning that the different possible behaviors correspond to
the three phases occurring in the mean field theory of directed
polymers with random weights in which the free energy is built from a
$b$-adic canonical complex cascade \cite{DES} [see Remark~\ref{commdiv}(4)].

Functional central limit theorems (Theorems \ref{div2} and \ref{div3})
add to the previous results on the almost sure behavior of $(F_n)_{n\ge1}$.
Let $W$ be real valued, and $(F_n)_{n\ge1}$
converge uniformly almost surely and in $L^2$ norm to a function $F$.
Then, under weak additional assumptions on the
moments of $W$ we prove that $(F_n-F)/\sqrt{\mathbb{E}((F_n-F)(1)^2)}$
converges in law as $n$ tends to $\infty$. Let $F_n$ be not bounded in $L^2$.
Then, under strong assumptions on the moments of $W$ we prove that
$F_n/\sqrt{\mathbb{E}(F_n(1)^2)}$ converges in law as $n$ tends to
$\infty$.

It is remarkable that in both cases
the weak limit is standard \cite{Mandfin}: it is Brownian motion $B$
in multifractal time $B\circ\widetilde F$ where
$\widetilde F$ is independent of $B$ and is the limit of a positive
canonical $b$-adic cascade. Thus the limit
process is one of the processes mentioned above as built from
statistically self-similar measures.
Among nontrivial functional central limit theorems, these results
seem to be the first to exhibit this kind of limit process. They reinforce
the importance of multifractal subordination in a monofractal process
as a natural operation.

The previous functional central limit theorems raise the following
question. Let $(F_n)_{n\ge1}$ strongly
converge to a nontrivial limit $F$. Can one---as it is the case for
some other classes of multifractal functions
\cite{Mandfin,JafMand2,BacryMuzy,Sadv}---decompose $F$ as $\widetilde B\circ\widetilde F$, where
$\widetilde B$ is a monofractal process and $\widetilde F$ the
indefinite integral of a statistically self-similar
measure? It turns out that there is a unique natural candidate
$\widetilde F$ as time change: There exists $\beta>1$ such that $\mathbb
{E}(|W|^\beta)= 1$ and $\widetilde F$ is the $b$-adic canonical cascade
generated by the martingale $(|Q_n|^\beta)_{n\ge1}$ . We give
sufficient conditions on the moments of $W$ for $F\circ\widetilde
F^{-1}$ to be indeed a monofractal (in the strong sense specified in
Section~\ref{1.1}) stochastic process of exponent $1/\beta$
(Theorem~\ref{timechange}). We do not know whether or not they are necessary.

We also consider a more general model of signed multiplicative
cascades, namely the signed extension of the most general
construction considered in \cite{M3}: We only assume that all the
vectors of the form $(W(w0),\dots,W(w(b-1)))$ are
independent and identically distributed, and $\mathbb{E}(\sum
_{i=0}^{b-1}W(i)/b)=1$. We call this class $b$-adic independent
cascades; such a cascade can be described as a $[0,1]$-martingale only
if $\mathbb{E}(W(i))>0$ for all $i$
(see Remark 2.3 in \cite{BJMpartI}). All our previous results have an
extension to this class. Moreover, one important additional case
appears in this class with respect to canonical cascades. It consists
of the conservative cascades for which $\sum_{i=0}^{b-1}W(i)/b=1$
almost surely. We can build examples of such conservative cascades
whose limit cannot be naturally decomposed as a monofractal function in
a multifractal time [see Theorems~\ref{conv-2}(2) and~\ref{notimechange}
and Remark~\ref{rem2.2}(5)].

\subsection{Organization of the paper}
The $b$-adic independent cascades are constructed in Section~\ref{CONST}.
The results regarding the uniform convergence and the representation of
the limit as a monofractal function in multifractal
time are given in Section~\ref{CONV} while results on degeneracy and
divergence as well as central limit theorems are
given in Section~\ref{DEGE}. Section~\ref{STABLE} provides a statement
on the multifractal nature of the limit whenever
it exists (proof and extensions are given in \cite{BJ}) as well as the
connection of the statistically self-similar processes
constructed in this paper with random variables and processes stable
under random weighted mean. Sections \ref{sec-5} and
\ref{secdiv} are dedicated to the proofs.

\subsection{Definitions}
Given an integer $b\geq2$, we denote by $\mathscr{A}$ the alphabet $\{
0,\dots,b-1\}$ and we define
$\mathscr{A}^*=\bigcup_{n\ge0} \mathscr{A}^n$ (by convention $\mathscr
{A}^0$ is the set reduced to the empty word denoted $\varnothing$). The
word obtained by concatenation of $u$ and $v$ in $\mathscr{A}^{*}$ is
denoted $u\cdot v$ and sometimes $uv$. For every $n\ge0$, the length
of an element of $\mathscr{A}^{n}$ is by definition equal to $n$ and
will be denoted by $|w|$. Let $n\ge1$ and $w=w_1\cdots w_n\in\mathscr
{A}^{n}$. Then for every $1\le k\le n$, the word $w_1\cdots w_k$ is
denoted $w|k$; if $k=0$, then $w|0$ stands for $\varnothing$.

For $w \in\mathscr{A}^{*}$, define $t_w=\sum_{i=1}^{|w|} w_i b^{-i}$
and $I_w=[t_w,t_w+b^{-|w|}]$.

If $f\in\mathcal{C}([0,1])$ and $I$ is a subinterval of $[0,1]$,
$\Delta f(I)$ denotes by the increment of $f$ over $I$. Also,
$\|f\|_\infty$ denotes the norm $\sup_{t\in[0,1]} |f(t)|$.

Denote by $(\Omega,\mathcal{B},\mathbb{P})$ the probability space on
which the random variables considered in this paper are
defined. Write $U\equiv V$ to express that the random variables $U$ and
$V$ have the same probability distribution.
The probability distribution of a random variable $V$ is denoted by
$\mathcal{L}(V)$.

\section{Construction of $b$-adic independent cascades and main
results}\label{badic}
\subsection{Construction}\label{CONST}

Let $W=(W_0,\dots,W_{b-1})$ be a complex vector whose components are
integrable and satisfy $\mathbb{E}(\sum_{i=0}^{b-1}W_i)=1$ (we have
modified the normalization with respect to the discussion of
Section~\ref{intro}). Then let $(W(w))_{w\in\mathscr{A}^{*}}$ be a
sequence of independent copies of $W$, and consider the sequence of
random functions
%
\begin{equation}\label{Fn'}
F_n(t)=F_{W,n}(t)=\int_0^t b^{n}\prod_{k=1}^{n}W_{u_{k}}(u|k-1) \,\mathrm{d}u
\end{equation}
(where each non $b$-adic point $u$ is identified with the word
$u_1\cdots u_n\cdots$ defined by the $b$-adic expansion $u=\sum_{k\ge
1}u_kb^{-k}$, and we recall that $u|k-1=u_1\cdots u_{k-1}$).
A special case playing an important role in the sequel is the
conservative one, that is, when $\sum_{i=0}^{b-1}W_i=1$ almost surely.
Remarkable functions obtained as the limit of such deterministic
sequences $F_n$ are the self-affine functions considered, for instance,
in \cite{Kono1,Bed,Urb1} (these functions are called self-affine
because their graphs are self-affine sets).

For $p\in\mathbb{R}_+$ let
%
\begin{equation}\label{varphi}
\varphi_W(p)=-\log_b\mathbb{E} \Biggl(\sum_{i=0}^{b-1}|W_i|^p \Biggr).
\end{equation}
Our assumptions imply $-1\le\varphi_W(0) \wedge\varphi_W(1) \le
\varphi_W(0) \vee\varphi_W(1) \le0$.

\begin{figure}[b]

\includegraphics{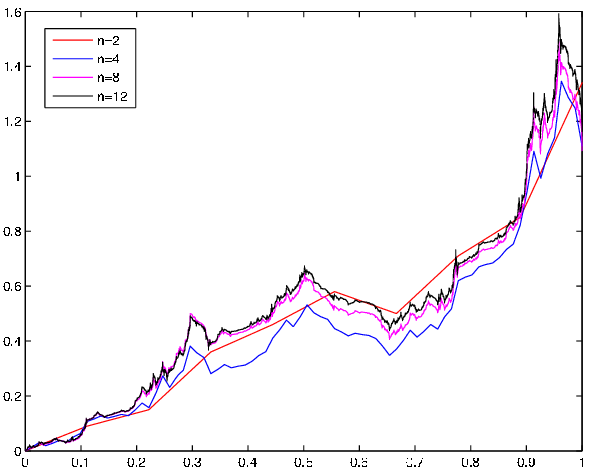}

\caption{Uniform convergence in the nonconservative case: $F_{W,n}$
for $n=2, 4, 8, 12$ in the case $b=3$ and $\varphi_W(\beta)=0$ for
$\beta\approx1.395$.}
\label{fig:1}
\end{figure}

\begin{figure}[t]

\includegraphics{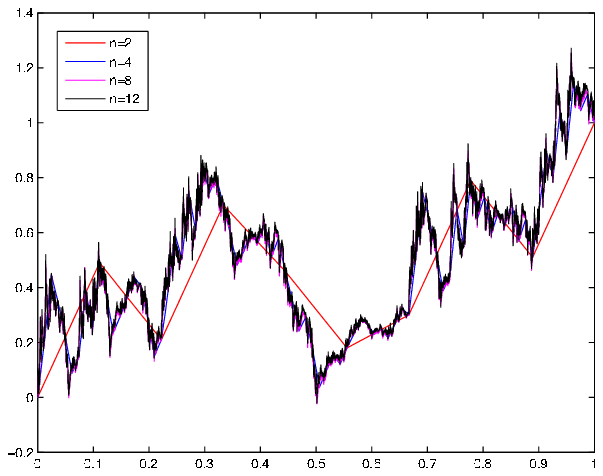}

\caption{Uniform convergence in the noncritical conservative case:
$F_{W,n}$ for $n=2, 4, 8, 12$ in the case $b=3$ and $\varphi_W(\beta
)=0$ for $\beta\approx2.172$. }
\label{fig:2}
\end{figure}
%
\begin{figure}[b]

\includegraphics{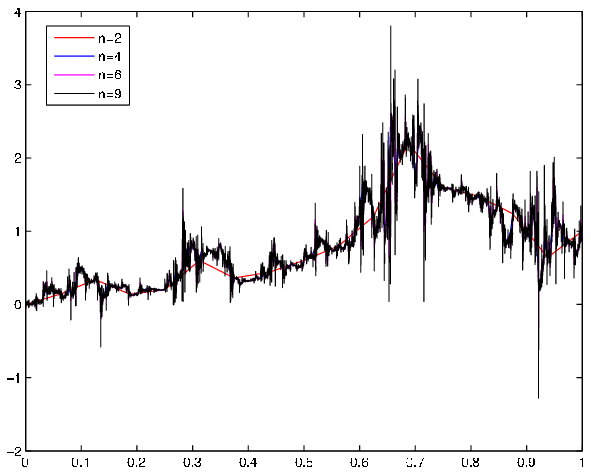}

\caption{Uniform convergence in the critical conservative case:
$F_{W,n}$ for $n=2, 4, 6, 9$ in the case $b=4$ and $\varphi_W(p)<0$
on $\mathbb{R}_+$, $\varphi_W(p)\to0$ when $p\to\infty$. The limit is
not uniformly H\"older.}
\label{fig:3}
\end{figure}

\subsection{Strong uniform convergence}\label{CONV}
We give sufficient conditions for the almost sure uniform convergence
of $(F_n)_{n\ge1}$. The asymptotic behavior of $(F_n)_{n\ge1}$ when
these conditions do not hold will be examined in Section~\ref{DEGE}.

We distinguish the conservative and nonconservative cases. Our results
are illustrated in Figures~\ref{fig:1} to \ref{fig:3}.

\begin{theorem}[(Nonconservative case)]\label{conv-1}
Suppose that $ \mathbb{P} (\sum_{i=0}^{b-1}W_i\neq1 )>0$
and there exists $p>1$ such that $\varphi_W(p)>0$. Suppose, moreover,
that either $p\in(1,2]$ or $\varphi_W(2)>0$.
\begin{longlist}[(1)]
\item[(1)]
$(F_n)_{n\ge1}$ converges uniformly, almost surely and in $L^p$ norm,
as $n$ tends to $\infty$, to a function $F=F_W$, which is nondecreasing if $W\ge0$. Moreover, the function $F$ is $\gamma$-H\"older
continuous for all $\gamma$ belonging to\break $(0,\max_{q\in(1,p]}\varphi_W(q)/q)$.

\item[(2)]$F$ satisfies the statistical scaling invariance property
%
\begin{equation}\label{sesi}
F=\sum_{i=0}^{b-1}\mathbf{1}_{[i/b,(i+1)/b]} \bigl(F(i/b)+W_i F_i\circ
S_i^{-1} \bigr),
\end{equation}
where $S_i(t)=(t+i)/b$, the random objects $W$, $F_0,\dots,F_{b-1}$ are
independent, and the $F_i$ are distributed like $F$ and the equality
holds almost surely.
\end{longlist}
\end{theorem}

\begin{theorem}[(Conservative case)]\label{conv-2} Suppose that
$ \mathbb{P} (\sum_{i=0}^{b-1}W_i=1 )=1$.
\begin{longlist}[(1)]
\item[(1)]
If there exists $p>1$ such that $\varphi_W(p)>0$, then the same
conclusions as in Theorem~\textup{\ref{conv-1}} hold.

\item[(2)] \textup{Critical case}. Suppose that $\lim_{p\to\infty} \varphi
_W(p)=0$ [in particular $\varphi_W$ is increasing and $\varphi_W(p)<0$
for all $p>1$]. This is equivalent to the fact that $\mathbb{P}(\forall
 0\le i\le b-1, |W_i|\le1)=1$ and $\sum_{i=0}^{b-1}\mathbb{P}(|W_i|=1)=1$.

Suppose also that $\mathbb{P}(\#\{i\dvtx |W_i|=1\}=1)<1$, and there exists
$\gamma\in(0,1)$ such that, with probability 1, one of the two
following properties holds for each $0\le i\le b-1$:
%
\begin{equation}\label{critical}
\cases{
\mbox{either } |W_i|\le\gamma,\cr
\mbox{or } |W_i|=1 \mbox{ and } (\sum_{k=0}^{i-1} W_i,\sum_{k=0}^i W_i
) \in\{(0,1),(1,0)\}.
}
\end{equation}
Then, with probability 1, $(F_n)_{n\ge1}$ converges almost surely
uniformly to a limit $F=F_W$ which is not uniformly H\"older and
satisfies part \textup{(2)} of Theorem~\ref{conv-1}.
\end{longlist}
\end{theorem}


%

\begin{remark}\label{Fn0}
(1) The sufficient condition for the convergence in $L^p$ of
complex-valued martingales like $(F_n(1))_{n\ge1}$ is known in the
context of martingales in the branching random walk (\cite
{Biggins,B3}); however, the sequence of functions $(F_n)_{n\ge1}$ is
not considered in these papers. When $W$ has nonnegative components,
it follows from \cite{DL} and \cite{KP} that this condition is necessary.

(2) Theorem~\ref{conv-2} goes beyond the construction of
deterministic self-affine functions (\cite{Urb1,Bed}) which all fall in
Theorem~\ref{conv-2}.

(3) The following discussion will be useful for the statement and proof
of Theorem~\ref{div}. It is easily seen that $F_n$ vanishes ($F_n=0$)
if and only if $\prod_{k=1}^{n}W_{w_{k}}(w|k-1)=0$ for all $w\in
\mathscr{A}^{n}$, and in this case, $F_k=0$ for all $k>n$. Thus, if we
denote by $\mathscr{V}$ the event $\{\exists n\ge1\dvtx  F_n=0\}$, we
have $\mathscr{V}=\liminf_{n\to\infty} \{F_n=0\}$. Notice that $\mathbb
{P}(\mathscr{V})=0$ in the conservative case.

By construction, there are $b$ independent copies $(F_{i,n})_{n\ge1}$
of $(F_n)_{n\ge1}$, independent of $W$, and converging respectively to
$F_i$ almost surely, such that for $n\ge1$ we can write
\[
F_n=\sum_{i=0}^{b-1}\mathbf{1}_{[i/b,(i+1)/b]} \bigl(F_n(i/b)+W_i
F_{i,n-1}\circ S_i^{-1} \bigr).
\]
Thus, $ \{F_n=0\}=\bigcap_i \{W_i F_{i,n-1}=0\}$. Similarly, $ \{F=0\}
=\bigcap_i \{W_i F_{i}=0\}$. It is then not difficult to see that under
the assumptions of Theorem \ref{conv-1}, $\mathbb{P}(\mathscr{V})$ and
$\mathbb{P}(F=0)$ are equal to the unique fixed point distinct from 1
of the convex polynomial
\[
P(x)=\sum_{k=0}^{b}\mathbb{P} (\#\{0\le i\le b-1\dvtx |W_i|>0\}=k) x^k.
\]
Consequently, since $\mathscr{V}\subset\{F=0\}$, these events differ
from a set of null probability.

We can also interpret the set of words $w\in\mathscr{A}^n$ such that
$\prod_{k=1}^{n}W_{w_{k}}(w|k-1)=0$ as the nodes of generation $n$ of a
Galton--Watson tree whose offspring distribution is given by that of
the integer $N=\#\{0\le i\le b-1\dvtx |W_i|>0\}$. Then $P(x)=\mathbb{E}(x^N)$.
\end{remark}

We end this section by introducing auxiliary nonnegative $b$-adic
independent cascades which will play an important role in the rest of
the paper.
\begin{definition}\label{Wbeta}
If $\beta>0$ and $\varphi_W(\beta)>-\infty$, then for $w\in\mathscr
{A}^{*}$ let
\[
W^{(\beta)}(w)=b^{\varphi_W(\beta)}(|W_0(w)|^\beta,\dots
,|W_{b-1}(w)|^\beta),
\]
and simply denote $W^{(\beta)}(\varnothing)$ by $W^{(\beta)}$. We have
$\varphi_{W^{(\beta)}}(p)=\varphi_W(\beta p)-p\varphi_W(\beta)$ for all
$p>0$. In particular, $\varphi_{W^{(\beta)}}(1)=0$. If $\varphi
_{W^{(\beta)}}(p)>0$ for some $p\in(1,2)$, we denote by $F_{W^{(\beta
)}}$ the nondecreasing function obtained in Theorem \ref{conv-1} as
the almost sure uniform limit of $F_{W^{(\beta)},n}\dvtx t\in[0,1]\mapsto
\int_0^t b^{n}\prod_{k=1}^{n}W^{(\beta)}_{u_{k}}(u|k-1)\, \mathrm{d}u$.
\end{definition}

\subsection{Representation as a monofractal function in multifractal time}
\begin{figure}[t]

\includegraphics{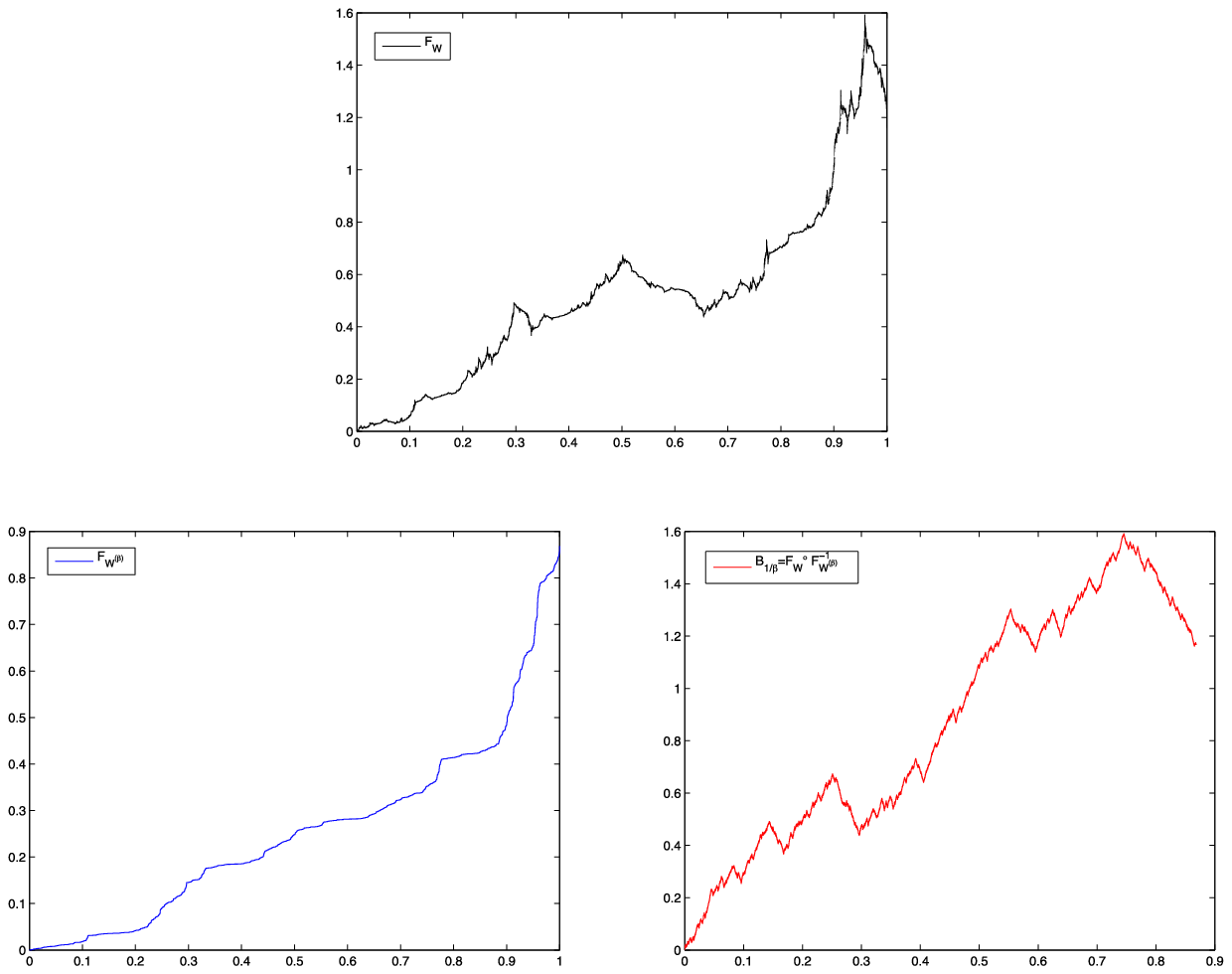}

\caption{The limit function $F_W$ corresponding to the construction of
Figure~\protect\ref{fig:1} (top) can be written as the monofractal function
$B_{1/\beta}=F_{W}\circ F_{W^{(\beta)}}^{-1}$ of exponent $H=1/\beta
\approx0.7168$ (right-bottom) in multifractal time $F_{W^{(\beta)}}$
(left-bottom).}
\label{fig:4}\vspace*{-1pt}
\end{figure}

\begin{figure}[t]

\includegraphics{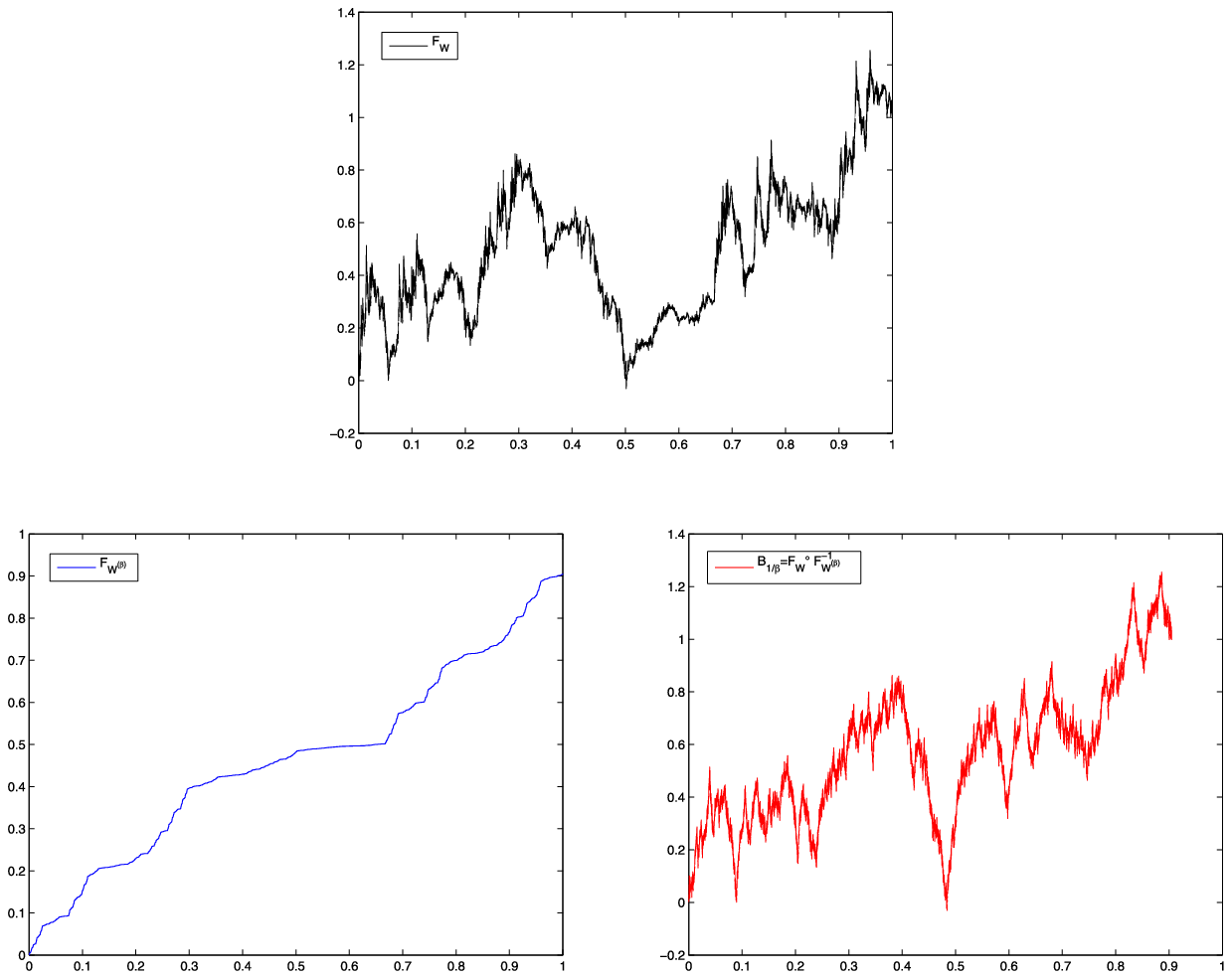}

\caption{The limit function $F_W$ corresponding to the construction of
Figure~\protect\ref{fig:2} (top) can be\vspace*{1pt} written as the monofractal function
$B_{1/\beta}=F_{W}\circ F_{W^{(\beta)}}^{-1}$ of exponent $H=1/\beta
\approx0.4604$ (right-bottom) in multifractal time $F_{W^{(\beta)}}$
(left-bottom).}\vspace*{-1pt}
\label{fig:5}
\end{figure}

As explained in the introduction, in order to qualitatively compare the
strong limit $F_W$ obtained in Theorem~\ref{conv-1} with other models
of multifractal processes, it is important to study the possibility to
decompose it as a monofractal function in multifractal time. Under the
assumptions of Theorems \ref{conv-1} and  \ref{conv-2}(1), if we
denote by $\beta$ the smallest solution of $\varphi_W(p)=0$, the only
natural choice at our disposal as time change is the function
$F_{W^{(\beta)}}$ introduced in Definition~\ref{Wbeta}. In the
deterministic case, it is elementary to check that $B_{1/\beta
}=F_W\circ F_{W^{(\beta)}}^{-1}$ is monofractal (Section~4.7 in \cite
{Mandfin}) in the strong sense. In the random case, this also true
under strong assumptions on the moments of $W$ as shows Theorem~\ref
{timechange} which is illustrated in Figures~\ref{fig:4} and \ref
{fig:5}. We do not know whether or not weaker assumptions on $W$ lead
to situations in which $B_{1/\beta}$ is not monofractal. However, the
functions constructed in Theorem~\ref{conv-2}(2) provide simple examples
of statistically self-similar continuous functions for which it seems
to be impossible to find a natural decomposition as monofractal
functions in multifractal time; at least such a time change cannot be
obtained as limit of a positive $b$-adic independent cascade.


\begin{theorem}\label{timechange}
Suppose that $\mathbb{P}(W\in\mathbb{C}^b\setminus\mathbb{R}_+^b)>0$
and $\sum_{i=0}^{b-1}\mathbb{E}(|W_i|^p)<\infty$ for all $p\in\mathbb{R}$.

Suppose also that the assumptions of Theorems \ref{conv-1} or
\textup{\ref{conv-2}(1)} hold. Let $\beta$ be the smallest solution of
$\varphi_W(p)=0$, and suppose that $\varphi_W(p)>0$ for all $p>\beta$.
We have $\beta>1$, $\varphi_{W^{(\beta)}}(1)=0$ and $\varphi_{W^{(\beta
)}}(p)>0$ for all $p>1$. Let $B_{1/\beta}=F_W\circ F_{W^{(\beta)}}^{-1}$.

With probability 1, the function $B_{1/\beta}$ is a monofractal
function in the strong sense, with constant pointwise H\"older exponent
$1/\beta$.
\end{theorem}

\begin{theorem}\label{notimechange}
Suppose that the assumptions of Theorem~\textup{\ref{conv-2}(2)} hold and the
components of $W$ do not vanish. With probability 1, if the function
$F_W$ can be decomposed as a monofractal function $B$ in a multifractal
time $G$, then the H\"older exponent of $B$ is 0, and $G$ must be such
that $
\dim ([0,1]\setminus\overline E_G(\infty) )=0$, where
\[
\overline E_G(\infty)= \biggl\{t\in[0,1]\dvtx \limsup_{r\to0^+}\frac
{\log(\operatorname{Osc}_{G}([t-r,t+r]))}{\log(r)}=\infty\biggr\}.
\]
\end{theorem}

\begin{remark}\label{rem2.2}
(1) Notice that due to Theorem \ref{conv-1}, $1/\beta$ must belong to
$(1/2,1)$ when $\mathbb{P}(\sum_{k=0}^{b-1}W_k=1)<1$.

(2) When $W$ is deterministic, we are necessarily in the conservative
case $\sum_{i=0}^{b-1}W_i=1$, and Theorem \ref{timechange} is well
known (see, for instance, Section 4.7 in \cite{Mandfin}). In this case,
it is also the simplest illustration of the general result obtained in
\cite{Sadv} regarding the representation of multifractal functions as
monofractal functions (in the strong sense) in multifractal time (see
also \cite{JafMand2} for another illustration of this concept).

(3) Under the assumptions of Theorem \ref{timechange}, it seems
possible to obtain the result by using the general approach developed
in \cite{Sadv}. However, this necessitates the use of the extension to
the present context of some sophisticated estimates developed for
positive cascades in \cite{BS1}. Thus we will give a short and
self-contained proof.

(4) The moments of $ \frac{\operatorname{Osc}_{F_W}([0,1])}{F_{W^{(\beta
)}}(1)^{1/\beta}}$ are all finite under the strong assumptions of
Theorem \ref{timechange}. Suppose that we have found sufficient
conditions under which there exists $q\neq0$ such that $
\mathbb{E} [ (\frac{\operatorname{Osc}_{F_W}([0,1])}{F_{W^{(\beta)}}(1)^{1/\beta}}) ^q
]=\infty$. Then we know how to prove\vspace*{1pt} that the function $B_{1/\beta}$ is
not monofractal if $q>0$ and not strongly monofractal if $q<0$. Thus
finding information on the $ \frac
{\operatorname{Osc}_{F_W}([0,1])}{F_{W^{(\beta)}}(1)^{1/\beta}}$ moments behavior
under weak assumptions on $W$ remains an important open question to
complete the description of $B_{1/\beta}$.

(5) As a consequence of Theorem~\ref{notimechange}, we get that the
time change $G$ cannot be equivalent to the limit of a positive
$b$-adic independent cascade obtained as in Theorems~\ref{conv-1} or~\ref
{conv-2}, in the sense that their derivative in the distribution sense
are equivalent positive measures. Indeed, in this case the analysis of
such a measure achieved in \cite{B1} or  \cite{KP}, implies that there
exists $D>0$ such that $ \lim_{n\to0^+} \frac{\log_b (\operatorname{Osc}_{G}
(I_n(t)))}{-n}=D$ at each point $t$ of a set $E$ of Hausdorff dimension~$D$.

Also, we can notice that if a time change $G$ as described in
Theorem~\ref{notimechange} does exist, from the multifractal analysis
point of view, the decomposition would not simplify the study of $F_W$
since it would necessitate having a very fine description of the
pointwise divergence of $ \frac{\log(\operatorname{Osc}_{G}([t-r,t+r]))}{\log
(r)}$ inside the set $\overline E_G(\infty)$.

An example of increasing function $G$ such that $\dim ([0,1]\setminus
\overline E_G(\infty) )=0$ is obtained as follows:
\[
G(t)=\int_0^t\sum_{n\ge1}\sum_{0\le k<b^n}b^{-n^2} \bigl(b^{n^3}\mathbf
{1}_{[kb^{-n}, k b^{-n}+b^{-n^3}]}(u) \bigr)\, \mathrm{d}u.
\]
We leave the reader check that $t$ is not a $b$-adic number and if
\[
\delta_t=\limsup_{n\to\infty}\sup_{0\le k<b^n} \frac{\log
(|t-kb^{-n}|)}{-n}<\infty,
\]
then $\lim_{r\to0^+}\frac{\log(\operatorname{Osc}_{G}([t-r,t+r]))}{\log
(r)}=\infty$. Moreover, it is clear that $\dim \{t\dvtx\delta_t=\break \infty\}=0$.
\end{remark}

\subsection{Degeneracy, divergence and weak uniform convergence}\label{DEGE}
Recall that $\varphi_W$ is concave, $\varphi
_W(0)<0$ and $\varphi_W(1)\le0$. Let us define
%
\begin{equation}\label{paramp0}
p_0=\sup\{p\dvtx \varphi_W'(p)\mbox{ exists and }\varphi_W'(p)p-\varphi
_W(p)>0\}.
\end{equation}
In order to simplify the next discussion, we assume that $\varphi
_W(p)>-\infty$ for all $p\ge0$.

Since $\varphi_W(0)<0$, we have $p_0>0$. Also, if $p_0<\infty$ then
$\varphi_W(p_0)=0$ if and only if $\varphi_W(p_0)=\varphi_W'(p_0)=0$.
If $p_0=\infty$, we define $\varphi_W(p_0)=\lim_{p\to\infty}\varphi_W(p)$.

In the previous section, we have dealt with the convergence of $F_n$ in
the following cases that we gather in the condition \textup{(C)}:

\begin{longlist}
\item[(C):] One of the following three cases arises:
\begin{longlist}[(1)]
\item[(1)] There exists $p\in(1,2]$ such that $\varphi_W(p)>0$;\vspace*{1pt}

\item[(2)] $\mathbb{P}(\sum_{k=0}^{b-1}W_k=1)=1$ and there exists $p>1$ such
that $\varphi_W(p)>0$;\vspace*{2pt}

\item[(3)] $\mathbb{P}(\sum_{k=0}^{b-1}W_k=1)=1$,\vspace*{1pt} $\mathbb{P}(\forall 0\le
i\le b-1,\ |W_i|\le1)=1$, $\sum_{i=0}^{b-1}\mathbb{P}(|W_i|=1)= 1$ and
$\mathbb{P}(\#\{i\dvtx |W_i|=1\}=1)<1$. Equivalently, $\lim_{p\to\infty}
\varphi_W(p)=0$ and $\mathbb{P}(\#\{i\dvtx |W_i|=1\}=1)<1$.
\end{longlist}
\end{longlist}

Suppose that \textup{(C)} does not hold. We cannot have simultaneously
$p_0\in(1,2]$ and $\varphi_W(p_0)>0$. Also, $\varphi_W(2)\le0$.
Moreover, if $\mathbb{P}(\sum_{k=0}^{b-1}W_k= 1)=1$, then $\varphi
_W(p_0)\le0$, and if $p_0=\infty$, then $\sum_{i=0}^{b-1}\mathbb
{P}(|W_i|=1)>1$ or $\mathbb{P}(\#\{i\dvtx |W_i|=1\}=1)=1$ (see also
Remark~\ref{commdiv}).

  The following results concern the asymptotic behavior of $F_n$
when \textup{(C)} does not hold. Before stating it, we recall the
discussion of Remark \ref{Fn0}(3).

\begin{theorem}[(Degeneracy and divergence)] \label{div}
Suppose that \textup{(C)} does not hold and $\varphi_W(p)>-\infty$ for all
$p\ge0$.
\begin{longlist}[(1)]
\item[(1)] Suppose that $p_0\in(0,1]$. Then, for all $\alpha\le\varphi
_W(p_0)/p_0$, $b^{n\alpha}F_n$ converges almost surely uniformly to 0,
and for all $\alpha> \varphi_W(p_0)/p_0$, $b^{n\alpha}F_n$ is unbounded
almost surely, conditionally on $\mathscr{V}^c$.

\item[(2)] Suppose that $p_0\in(1,2]$. We have $\varphi_W(p_0)\le0$. Then,
for all $\alpha> \varphi_W(p_0)/p_0$, $b^{n\alpha}F_n$ is unbounded
almost surely, conditionally on $\mathscr{V}^c$. In particular, if
$\varphi_W(p_0)<0$, then $F_n$ is unbounded almost surely,
conditionally on $\mathscr{V}^c$.

\item[(3)] Suppose that $p_0> 2$ and $\mathbb{P}(\sum_{k=0}^{b-1}W_k= 1)=1$.
We have $\varphi_W(p_0)\le0$. If $p_0<\infty$, the same conclusions as
in \textup{(2)} hold. If $p_0=\infty$, then $(F_n)_{n\ge1}$ diverges in $\mathcal
{C}([0,1])$ almost surely.

Moreover, in both cases there is no sequence $(r_n)_{n\ge1}$ tending
to 0 or $\infty$, as $n\to\infty$, such that $r_n F_n$ converges in law
to a nontrivial limit in $\mathcal{C}([0,1])$.

\item[(4)] Suppose that $p_0> 2$ and $\mathbb{P}(\sum_{k=0}^{b-1}W_k\neq
1)>0$. For $n\ge1$, let $r_n=b^{n\varphi_W(2)/2}$ if $\varphi_W(2)<0$
and $r_n=n^{-1/2}$ if $\varphi_W(2)=0$. The probability distributions
of the random functions $r_n F_n$ form a tight sequence, and for all
$\alpha>\varphi_W(2)/2$, $b^{n\alpha}F_n$ is unbounded almost surely,
conditionally on $\mathscr{V}^c$.

\end{longlist}
\end{theorem}

\begin{remark}\label{commdiv}
(1) Suppose that \textup{(C)} does not hold and $p_0<\infty$. Thus $p_0\in
(1,2]$ and $\varphi_W(p_0)=0$, or when $\mathbb{P}(\sum_{k=0}^{b-1}W_k=
1)=1$, $p_0\in(1,\infty)$ and $\varphi_W(p_0)=0$. What we can only
prove is that $\lim_{n\to\infty}\sup_{w\in\mathscr{A}^{n}}|\Delta
F_n(I_w)|=0$ and $ \liminf_{n\to\infty} \frac{\log_b \sup
_{w\in\mathscr{A}^{n}}|\Delta F_n(I_w)|=0}{-n}$. This is not enough to
decide whether or not $F_n$ is convergent. It is mainly for the same
reason that we cannot deal with the case $\alpha= \varphi_W(p_0)/p_0$
in Theorem \ref{div}(2) and (3) (when $p_0<\infty$).

(2) When $p_0=\infty$, Theorem~\ref{div} tells nothing about the case
where the assumptions of Theorem~\ref{conv-1}(2) hold except that there
is no $\gamma\in(0,1)$ such that \eqref{critical} holds.

(3) The results obtained in \cite{DL,Gu} when $W\ge0$ show that in
this case, when $F_n$ converges almost surely uniformly to 0, there
does not exist a sequence $(a_n)_{n\ge1}$ such that $F_n/a_n$ converge
in law to a nontrivial process as $n$ tends to $\infty$ (see the
discussion in Section VIII of \cite{Gu}). Theorem~\ref{div2} shows that
allowing the components of $W$ to take values in $\mathbb{R}_-^*$
yields a completely different situation.

 (4) One referee brought to our attention reference \cite
{DES} which studies the possible phases in mean field theory of
directed polymers with random complex weights. Specifically, the
mathematical question discussed in \cite{DES} concerns the asymptotic
behavior of $\log|F_n(1)|/n$ in the canonical case, without special
assumption on the value of $\mathbb{E}(W_0)$, but under the assumption
that the distribution of $W_0$ is continuous, all its positive moments
are finite and $W_0/|W_0|$ is independent of $|W_0|$. The same
parameter $p_0$ as the one defined in \eqref{paramp0} is considered and
three phases are distinguished: phase I corresponds to \textup{(C)}, phase
II to parts (1) and (2) of Theorem~\ref{div} and phase III to part (4) of
Theorem~\ref{div}. The results obtained in \cite{DES} complete those
obtained in the present paper for phases II and III by showing that
$\log|F_n(1)|/n$ converges in probability to $-\varphi_W(p_0)/p_0$ in
phase II and $-\varphi_W(2)/2$ in phase III. In phase III and under the
assumptions of Theorem~\ref{div2}, our functional central limit theorem
describes the asymptotic behavior of the distribution of $\log|F_n(1)|/n$.
\end{remark}

Part (4) of Theorem~\ref{div} is restated and refined in Theorems~\ref
{div-caspart} and~\ref{div2}.

We now cease to assume that $\varphi_W(p)>-\infty$ for all $p\ge0$,
but assume that $\varphi_W(2)>-\infty$.

Define
\[
\sigma=
\cases{
 \sqrt{\dfrac{\mathbb{E}( |\sum_{i=0}^{b-1}W_i |^2
)-1}{\mathbb{E} (\sum_{i=0}^{b-1}|W_i|^2 )-1}},&\quad\mbox{if }$\varphi
_W(2)<0,$\vspace*{2pt}\cr
\sqrt{\sum_{i\neq j}\mathbb{E}(W_i\overline{W_j})},&\quad\mbox{if
}$\varphi_W(2)=0$.
}
\]

\begin{theorem}[{[Tightness of $(\mathcal{L}(F_n/ \sqrt{\mathbb{E}(F_n(1)^2)})_{n\ge1})$]}]\label{div-caspart}
Suppose that\break $\mathbb{P}(\sum_{k=0}^{b-1}W_k\neq1)>0$ and \textup{(C)}
does not hold. In particular, $\varphi_W(2)\le0$.
\begin{longlist}[(1)]
\item[(1)] The sequence $(F_n(1))_{n\ge1}$ is unbounded in $L^2$ norm.
Specifically, we\vspace*{1pt} have $\mathbb{E}(|F_n(1)|^2)\sim\sigma^2 b^{-n\varphi
_W(2)}$ if $\varphi_W(2)<0$ and $\mathbb{E}(|F_n(1)|^2)\sim\sigma^2 n$
if $\varphi_W(2)=0$.

\item[(2)] Suppose that $p_0>2$. Equivalently, $\varphi_W(p)/p>\varphi
_W(2)/2$ near $2^+$.\vspace*{1pt} For $n\ge1$ let $Z_n=F_n/ \sqrt{\mathbb{E}(F_n(1)^2)}$.

The sequence $(Z_n(1))_{n\ge1}$ is bounded in $L^p$ norm for all $p$
such that\break $\varphi_W(p)/p>\varphi_W(2)/2$.

Moreover, the probability distributions of the random continuous\vspace*{1pt}
functions $Z_n=F_n/ \sqrt{\mathbb{E}(F_n(1)^2)}$ form a tight sequence.

\item[(3)] Suppose that $p_0>2$. We have $\varphi_{W^{(2)}}(p)>0$ near $1^+$
(remember Definition~\ref{Wbeta}). Suppose, moreover, that $W$ is
$\mathbb{R}^b$-valued and $(Z_n)_{n\ge1}$ converges in law, as $n$
tends to $\infty$. Then the weak limit of $Z_n$ is the Brownian motion
in multifractal time $Z=B\circ F_{W^{(2)}}$ where $B$ is a standard
Brownian motion independent of $F_{W^{(2)}}$. Moreover, $Z$ satisfies
the statistical scaling invariance property,
%
\begin{equation}\label{funcBrF}
Z\equiv\sum_{i=0}^{b-1}\mathbf{1}_{[i/b,(i+1)/b]} \bigl(Z(i/b)+b^{\varphi
_W(2)/2}W_i Z_i\circ S_i^{-1} \bigr),
\end{equation}
where $S_i(t)=(t+i)/b$, the random objects $W$, $Z_0,\dots,Z_{b-1}$ are
independent, and the $Z_i$ are distributed like $Z$.
\end{longlist}
\end{theorem}

\begin{figure}[t]

\includegraphics{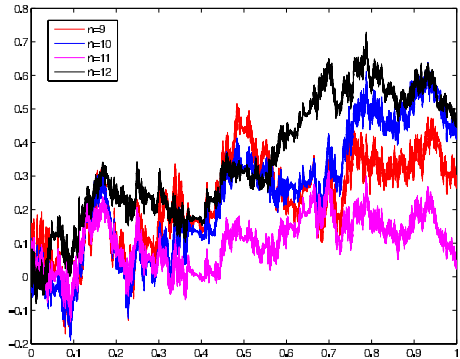}

\caption{Weak convergence of $F_{W,n}/(\sigma b^{-n\varphi_W(2)/2})$ to
a Brownian motion in multifractal time when $\varphi_W(2)\le0$ and
$\varphi_W(p)/p>\varphi_W(2)/2$ for all $p>2$ in the nonconservative case.}
\label{fig:6}\vspace*{-2pt}
\end{figure}

The following result completes Theorem~\ref{div-caspart} and is
illustrated in Figure~\ref{fig:6}.
\begin{theorem}[(Functional central limit theorem when
$F_n$ is unbounded)]\label{div2} Suppose that $\mathbb{P}(\sum_{k=0}^{b-1}W_k\neq
1)>0$ and \textup{(C)} does not hold. Suppose, moreover, that $W$ is
$\mathbb{R}^b$-valued and $\varphi_W(p)/p>\varphi_W(2)/2$ over
$(2,\infty)$ (equivalently $\varphi_W$ is increasing, or $|W_k|\le1$
almost surely for all $0\le k\le b-1$).

Then $(Z_n)_{n\ge1}$ converges in law, as $n$ tends to $\infty$, to
the Brownian motion in multifractal time $Z$ described in Theorem \ref
{div-caspart}\textup{(3)}. Also, the probability distribution of $Z(1)$ is
determined by its moments.
\end{theorem}

Theorem \ref{div2} has the following natural counterpart when
$(F_n)_{n\ge1}$ converges almost surely.
\begin{theorem}[(Functional central limit theorem when $F_n$ converges)]\label{div3}
Suppose that $\mathbb{P}(\sum_{k=0}^{b-1}W_k\neq
1)>0$ and \textup{(C)} holds. Suppose, moreover, that $W$ is $\mathbb
{R}^b$-valued and $\varphi_W(p)/p>\varphi_W(2)/2$ near $2^+$.

Then $(F_n-F)/\sqrt{\mathbb{E}((F_n-F)(1)^2)}$ converges in law, as $n$
tends to $\infty$, to the Brownian motion in multifractal time
described in Theorem \textup{\ref{div-caspart}(3)}. Moreover, $\sqrt{\mathbb
{E}((F_n-F)(1)^2)}=b^{-n\varphi_W(2)/2}\sqrt{\mathbb{E} ((1-F_W(1))^2 )}$.
\end{theorem}

\begin{remark}(1) The equivalent of $\mathbb{E}(|F_n(1)|^2)$ given in
Theorem \ref{div-caspart}(1) is valid in general when $\varphi_W(2)\le0$.

(2) In Theorem \ref{div2}, if the components of $W$ have the
same constant modulus, then the limit process is the standard Brownian
motion over $[0,1]$. In particular, we recover the result obtained in
\cite{BM}, where the components are i.i.d.

(3) When the components of $W$ are nonnegative, independent,
and identically distributed, the scalar central limit theorem deduced
for the term $(F_n-F)(1)/\sqrt{\mathbb{E}((F_n-F)(1)^2)}$ from Theorem
\ref{div2} is almost a restatement of the central limit theorem
established in Corollary 4.3 in \cite{OW}.\vspace*{-1pt}
\end{remark}


%

\subsection{Multifractal nature and scaling invariance properties}\label{STABLE}\vspace*{-1pt}

\subsubsection{Multifractal nature}\label{monomulti}
The multifractal analysis of the sample paths of the limit process $F$
obtained in Theorem~\ref{conv-1} is achieved in \cite{BJ}. Here we
state a simplified version of the result obtained in \cite{BJ}.
\begin{theorem}
Under the assumptions of Theorem~\ref{conv-1}, assume that\break
$\mathbb{E}(\sum_{i=0}^{b-1}|W_i|^q )$ is finite for all
$q\in\mathbb{R}$ and $W\neq(1/b,\dots,1/b)$ (otherwise $F$ is the
identity function of
$[0,1]$). Then, with probability 1, for all $h\ge0$, $\dim E_F(h)=\inf
_{q\in\mathbb{R}}h
q+\log_b\mathbb{E} (\sum_{i=0}^{b-1}|W_i|^q )=\tau_F^*(h)$, a negative dimension
meaning that the corresponding iso-H\"older set is empty. Moreover, the
function $F$ is
monofractal if and only if there exists $H\in(0,1)$ such that
$|W_i|=b^{-H}$ almost
surely for all $\le i\le b-1$. In this case, the pointwise H\"older
exponent is everywhere
equal to $H$.
\end{theorem}

\subsubsection{Link with processes stable under random weighted
mean} The scaling properties of the stochastic processes obtained in
Theorems~\ref{conv-1} and~\ref{div2} are reminiscent of the fact that
the distribution of their increments between 0 and 1 is invariant under
random weighted mean (using the terminology of \cite{M3}), or
equivalently that this distribution is a fixed point of a smoothing
transformation (using the terminology of \cite{DL}). Specifically, this
increment $\widetilde Z$ satisfies a functional equation of the form
%
\begin{equation}\label{mandeq}
\widetilde Z\equiv\sum_{i=0}^{b-1} \widetilde W_i\widetilde Z_i,
\end{equation}
where $\widetilde Z\equiv\widetilde Z_i$ for all $i$, and the random
variables $\widetilde W=(\widetilde W_0,\dots,\widetilde W_{b-1})$,
$\widetilde Z_0, \ldots, \widetilde Z_{b-1}$ are independent.

The square integrable solutions of (\ref{mandeq}) have been studied in
\cite{Rosler}. Thanks to \cite{Rosler}, we can conclude that: (1) When
$p\ge2$, the limit process $F_W$ obtained in Theorem~\ref{conv-1} is
the unique square integrable continuous stochastic process with
expectation the identity function of $[0,1]$ and satisfying \eqref
{sesi}; (2) the limit process $B\circ F_{W^{(2)}}$ obtained in
Theorem~\ref{div2} is the unique centered square integrable continuous
stochastic process satisfying (\ref{funcBrF}).

It is interesting to compare the structure of the process constructed
in Theorem~\ref{conv-1} with other processes naturally associated with
(\ref{mandeq}), namely symmetric L\'evy processes in multifractal time
which also satisfy (\ref{sesi}) and are obtained as follows: If
$\widetilde W\ge0$ and there exists $\beta\in(1,2)$ such that $\varphi
_{\widetilde W}(\beta)=0$ and $\varphi_{\widetilde W}'(\beta)>0$, then
it is noticed in \cite{DL,Gu,Liu2001} that a solution of (\ref{mandeq})
is $X_\beta\circ F_{\widetilde W^{(\beta)}}(1)$ where $F_{\widetilde
W^{(\beta)}}$ is the nondecreasing function constructed in
Definition~\ref{Wbeta}, and $X_\beta$ is a symmetric stable L\'evy
process of index $\beta$ independent of $F_{\widetilde W^{(\beta)}}$.
The multifractal nature of the jump processes $X_\beta\circ
F_{\widetilde W^{(\beta)}}$ is given in \cite{BS}.

\section[Proofs of Theorems 2.1, 2.2 and 2.3]{Proofs of Theorems~\protect\ref{conv-1},
\protect\ref{conv-2} and \protect\ref{timechange}}\label{sec-5}

We start with a remark. Under the assumptions of Theorems~\ref{conv-1}
and~\ref{conv-2}(1), we have $-1\le\varphi_W(0)<\varphi_W(1)\le
0<\varphi_W(p)$. Since $\varphi_W$ is concave this implies that $\varphi
_W(q)<q-1$ for all $q\in(1,p]$ except if $\varphi_W(q)=q-1$ for all
$q$. This can happen only if the components of $W$ are positive and
equal to $1/b$ almost surely. In this case $F_n(t)=t$ for all $n\ge1$
and $t\in[0,1]$ and the result obviously holds. We exclude this case
in the rest of this section.

For $w\in\mathscr{A}^{*}$, we denote by $(F^{[w]}_n)_{n\ge1}$ the copy
of $(F_n)_{n\ge1}$ constructed with the random vectors $(W(w\cdot
u))_{u\in\mathscr{A}^{*}}$:
\[
F^{[w]}_n(t)=\int_0^t b^{n}\prod_{k=1}^{n}W_{u_{k}}(w\cdot u||w|+k-1)
\,\mathrm{d}u.
\]
For $n> |w|$, the increment $\Delta F_n(I_w)$ of $F_n$ over $I_w$ takes the~form
%
\begin{equation}\label{self-sim'}
\Delta F_n(I_w)=Q(w) F_{n-|w|}^{[w]}(1),
\end{equation}
where
\[
Q(w)=\prod_{k=0}^{|w|-1}W_{w_{k+1}}(w|k).
\]
This implies in particular that for every $n\ge1$ we have
%
\begin{equation}\label{foncmand}
F_n=\sum_{i=0}^{b-1}\mathbf{1}_{[i/b,(i+1)/b]}
\bigl(F_n(i/b)+W_iF_{n-1}^{[i]}\circ S_i^{-1} \bigr).
\end{equation}
Moreover, $Q(w)$ and $F_{n-|w|}^{[w]}(1)$ are independent, and for each
$p\ge1$, the families $\{ F^{[w]}_n\}_{n\ge1}$, $w\in\mathscr{A}^p$,
are independent.

\subsection[Proofs of Theorems 2.1 and 2.2(1)]{Proofs of Theorems~\protect\ref{conv-1}
and \protect\ref{conv-2}(1)}
If $\mathbb{P}(\sum_{i=0}^{b-1}W_i=1)=1$ then $F_k(1)=1$ almost surely.
If $p\in(1,2]$, the fact that the martingale $(F_k(1))_{k\ge1}$
converges almost surely and in $L^p$ norm is a consequence of Theorem 1
in \cite{Biggins}, and the case $p>2$ is a consequence of Theorem 1 in
\cite{B3} (the positive case is treated in \cite{DL} and \cite{KP}).
Then, equation \eqref{self-sim'} implies the almost sure convergence of
the $b$-adic increments of $F_n$.


Now we establish the almost sure uniform convergence of $F_n$. When
$F_n$ can be interpreted as a $[0,1]$-martingale, that is when the
components of $W$ have positive expectations (see \cite{BJMpartI}), the
proof provides a simpler alternative to the general proof given in \cite
{BJMpartI}.

Let $q\in(1,p]$ such that $\varphi_W(q)>0$ and define $M_q=\mathbb
{E}(\sup_{k\ge1}|F_k(1)|^q)$. By using (\ref{self-sim'}) as well as
the martingale property of $(F_k(1))_{k\ge1}$ and Doob's inequality we get
\begin{eqnarray*}
\mathbb{E}\Bigl(\sup_{n\ge1}|\Delta F_n(I_w)|^q\Bigr)&\le &\sum_{n=1}^{|w|}\mathbb
{E}(|\Delta F_n(I_w)|^q)+\mathbb{E}(|Q(w)|^q)\mathbb{E}\Bigl(\sup
_{n>|w|}\bigl|F_{n-|w|}^{[w]}(1)\bigr|^q\Bigr)
\\
&\le& \sum_{n=1}^{|w|} b^{-(|w|-n)q}\mathbb{E}(|Q(w|n)|^q)+
C_qM_q\mathbb{E}(|Q(w)|^q)
\end{eqnarray*}
for some constant $C_q$. Consequently, for $\gamma>0$ and $N\ge1$ we have
\begin{eqnarray*}
&&\mathbb{P}\Bigl(\max_{w\in\mathscr{A}^{N}}\sup_{n\ge1}|\Delta
F_n(I_w)|>b^{-\gamma N}\Bigr)
\\
&&\qquad\le b^{\gamma Nq}\sum_{w\in\mathscr{A}^{N}}\sum_{n=1}^N
b^{-(N-n)q}\mathbb{E}(|Q(w|n)|^q)+ C_qM_q \mathbb{E}(|Q(w)|^q)
\\
&&\qquad= b^{\gamma Nq} \Biggl[\sum_{n=1}^N b^{-(N-n)(q-1)}b^{-n\varphi
_W(q)}+C_qM_q b^{-N\varphi_W(q)} \Biggr]
\\
&&\qquad \le \biggl[ \frac{b^{q-1-\varphi_W(q)}}{b^{q-1-\varphi_W(q)}-1} +C_q M_q \biggr]
b^{\gamma Nq}b^{-N\varphi(q)},
\end{eqnarray*}
where we used the fact that $\varphi_W(q)< q-1$. It follows that if
$\gamma<\varphi_W(q)/q$, we have
\[
\sum_{N\ge1}\mathbb{P}\Bigl(\max_{w\in\mathscr{A}^{N}}\sup_{n\ge1}|\Delta
F_n(I_w)|>b^{-\gamma N}\Bigr)<\infty.
\]
Due to the Borel--Cantelli lemma, we conclude that, with probability 1,
%
\begin{equation}\label{ascoli}
\mbox{for $N$ large enough,$\qquad \max\limits_{w\in\mathscr{A}^{N}}\sup\limits_{n\ge
1}|\Delta F_n(I_w)|\le b^{-\gamma N}$}.
\end{equation}
Next, we use the following classical property: for any continuous
complex function~$f$ on~$[0,1]$, one has
%
\begin{equation}\label{holder'}
\omega(f,\delta) \le2b \sum_{n\ge-{\log\delta/\log b}}
\sup_{w\in\mathscr{A}^{n}} \Delta f (I_w),
\end{equation}
where $\omega(f,\delta)$ stands for the modulus of continuity of $f$,
\begin{eqnarray*}
\omega(f,\delta)=\mathop{\sup_{t,s\in[0,1]}}_{ |t-s|\le
\delta}|f(t)-f(s)|.
\end{eqnarray*}
Since $F_n(0)=0$ almost surely for all $n\ge1$, it follows from \eqref
{ascoli}, \eqref{holder'} and Ascoli--Arzela's theorem that, with
probability 1, the sequence of continuous functions $(F_n)_{n\ge1}$ is
relatively compact, and all the limit of subsequences of $F_n$ are
$\gamma$-H\"older continuous for all $\gamma<\max_{q\in(1,p]}\varphi
(q)/q$. Moreover, by the self-similarity of the construction \eqref
{self-sim'} we know that $F_n$ converges almost surely on set of
$b$-adic points. This yields the uniform convergence of $F_n$ and the
H\"older regularity property of the limit~$F$.

To see that $(F_n)_{n\ge1}$ converges in $L^p$ norm, it is enough to
prove that the sequence $ (\mathbb{E}(\sup_{1\le k\le n}\|F_k\|_\infty
^p))_{n\ge1}$ is bounded.

For $n\ge1$ and $0\le i\le b-1$ define
\begin{eqnarray*}
S_n=\sup_{1\le k\le n}\|F_k\|_\infty, \qquad  S_{n}(i)=\sup_{1\le k\le n}\bigl\|
F_k^{[i]}\bigr\|_{\infty}
\end{eqnarray*}
 and
\begin{eqnarray*}
  \widetilde S_n(i)=\sup_{1\le k\le
n}|F_k(ib^{-1})|.
\end{eqnarray*}
Due to (\ref{foncmand}) we have $ S_n\le\max_{0\le i\le b-1}
[\widetilde S_n(i) +|W_i| S_{n-1}(i) ]$, so
\[
\mathbb{E}( S_n^p )\le\sum_{i=0}^{b-1}\mathbb{E}
\bigl([|W_i|S_{n-1}(i)+\widetilde S_n(i)]^p \bigr).
\]

Denote by $\bar{p}\geq2$ the unique integer such that $\bar{p}-1<p\leq
\bar{p}$. By using the subadditivity of the mapping $x\ge0\mapsto x^{
p/\bar p}$ we get
\begin{eqnarray*}
&&\mathbb{E}\bigl(\bigl(|W_i|S_{n-1}(i)+\widetilde S_n(i)\bigr)^p\bigr)
 \\
 &&\qquad \leq\mathbb{E}
\bigl(\bigl((|W_i|S_{n-1}(i))^{p/\bar p}+\widetilde S_n(i)^{p/\bar p}\bigr)^{\bar p} \bigr)
\\
&&\qquad \le\mathbb{E}( |W_i|^p)\cdot \mathbb{E}(S_{n-1}(i)^p) + \mathbb
{E}(\widetilde S_n(i)^p)
\\
&&\qquad\quad {}+\sum_{m=1}^{\bar p-1} \left({\bar p\atop m}\right)\mathbb{E}\bigl(
[|W_i|S_{n-1}(i)]^{mp/\bar p} [\widetilde S_{n-1}(i)]^{(\bar p-m)
p/\bar p} \bigr).
\end{eqnarray*}
Now let us make some remarks:

\begin{longlist}[--]
\item[--] The H\"older inequality yields for any pair of nonnegative random
variables $(U,V)$ and $m\in[1,\bar p-1]$
\[
\mathbb{E}\bigl(U^{mp/\bar p}V^{(\bar p-m) p/\bar p}\bigr)\le\mathbb
{E}(U^p)^{m/\bar p} \mathbb{E}(V^p)^{(\bar p-m)/\bar p}.
\]

\item[--] The convergence in $L^p$ norm of $F_n(1)$ implies that $(\widetilde
S_n(i))_{n\ge1}$ is bounded in~$L^p$.

\item[--] The random variables $|W_i|$ and $S_{n-1}(i)$ are independent and
$|W_i|$ is in $L^p$.

\item[--] Since the expectation of $F_k(1)$ is equal to 1 for all $k\ge1$, we
have $1\le\mathbb{E}(S_{n-1}^p)^{m/\bar p}\le\mathbb
{E}(S_{n-1}^p)^{(\bar p-1)/\bar p}$ for all $m\in[1,\bar p-1]$.
\end{longlist}

The previous remarks imply the existence of two constants $A$ and $B$
independent of $i$ such that
\begin{eqnarray*}
&& \mathbb{E}\bigl(\bigl(W_iS_{n-1}(i)+\widetilde S_n(i)\bigr)^p\bigr)
\\
&&\qquad \leq \mathbb{E}( |W_i|^p) \mathbb{E}(S_{n-1}^p) +B+\sum_{m=1}^{\bar
p-1} \left({\bar p\atop m}\right) A^{m/\bar p}B^{(\bar p-m)/\bar p} \mathbb
{E}(S_{n-1}^p)^{(\bar p-1)/\bar p}
\\
&&\qquad = \mathbb{E}( |W_i|^p) \mathbb{E}(S_{n-1}^p) +B+(A+B)^{\bar p}\mathbb
{E}(S_{n-1}^p)^{(\bar p-1)/\bar p}.
\end{eqnarray*}
Summing this inequality over $i$ we find
\begin{eqnarray*}
\mathbb{E}( S_n^p ) \leq b^{-\varphi_W(p)} \mathbb{E}( S_{n-1}^p ) +
b(A+B)^{\bar p}\mathbb{E}(S_{n-1}^p)^{(\bar p-1)/\bar p}+ b\cdot B.
\end{eqnarray*}
For $x\ge0$ let $f(x)=b^{-\varphi_W(p)} \cdot x+ b(A+B)^{\bar
p}x^{(\bar p-1)/\bar p}+b\cdot B$. Since $b^{-\varphi_W(p)}<1$ and
$(\bar p-1)/\bar p<1$, there exists $x_0$ such that $f(x)<x$ for any
$x>x_0$,\vspace*{1pt} which implies $\mathbb{E}( S_n^p ) \leq\max\{x_0,\mathbb{E}(
S_{n-1}^p )\}$. This yields the conclusion.

Property \eqref{sesi} is a consequence of (\ref{foncmand}).

\subsection[Proof of Theorem 2.2(2)]{Proof of Theorem~\protect\ref{conv-2}(2)} By construction, if there
exists $0\le i\le b-1$ such that $\mathbb{P}(|W_i|>1)>0$, then $\lim
_{p\to\infty}\varphi_W(p)=-\infty$. Otherwise, we have $\lim_{p\to\infty
}\varphi_W(p)=-\log_b\sum_{i=0}^{b-1}\mathbb{P}(|W_i|=1)$, and by
concavity of $\varphi_W$ we have $\varphi_W(p)<0$ for all $p>0$ and
$\lim_{p\to\infty}\varphi_W(p)=0$ if and only if $\sum
_{i=0}^{b-1}\mathbb{P}(|W_i|=1)=1$.

For $n\ge1$ let us define $m_n=\max_{w\in\mathscr{A}^n}|Q(w)|$. Then,
due to \eqref{self-sim'}, for $p\ge1$, we have
\[
\|F_{n+p}-F_n\|_\infty\le\sup_{w\in\mathscr{A}^n}
|Q(w)| \bigl\|F^{[w]}_{p}-1\bigr\|_\infty\le m_n \sup_{w\in\mathscr{A}^n}\bigl(1+\bigl\|
F^{[w]}_{p}\bigr\|_\infty\bigr).
\]

\noindent
We are going to prove that $\lim_{n\to\infty} m_n=0$, and there exists
$C>0$ such that $\|F_p\|_\infty\le C$ almost surely for all $p\ge1$.
Thus $(F_n)_{n\ge1}$ is almost surely a Cauchy sequence.

We start with the proof of $\lim_{n\to\infty} m_n=0$. Due to the fact
that the components of $W$ are either equal to 1 or less than or equal
to $\gamma$, the sequence $(m_n)_{n\ge1}$ is nonincreasing and
$m_{n+1}=m_n$ or $m_{n+1}\le\gamma m_n$. Thus if $\lim_{n\to\infty}
m_n>0$, we can find an infinite word $w_1\cdots w_n \cdots$ in $
\mathscr{A}^{\mathbb{N}_+}$ and $n_0\ge1$ such that or all $n\ge n_0$,
$W_{w_{n+1}}(w_1\cdots w_n)=1$. By construction, such an infinite word
must belong to the boundary of a Galton--Watson tree rooted at
$w_1\cdots w_{n_0}$, and whose offspring distribution generating
function is given by $x\mapsto\sum_{k=0}^b\mathbb{P}(\#\{i\dvtx  |W_i|=1\}
=k) x^k$. This tree is subcritical, since we assumed that $\sum
_{k=1}^bk\mathbb{P}(\#\{i\dvtx |W_i|=1\}=k) =\sum_{i=0}^{b-1}\mathbb
{P}(|W_i|=1)=1$ and $\mathbb{P}(\#\{i\dvtx |W_i|=1\}=1)<1$. Consequently
its boundary is empty almost surely, and $m_n$ tends to 0 as $n\to\infty$.

Now we prove by induction that, with probability 1, for all $p\ge1$
and $w\in\mathscr{A}^p$, we have
%
\begin{equation}\label{fp}
\quad\quad \max\bigl(|F_p(t_{w})|, |F_p(t_w+b^{-p})| \bigr)\le1+\frac{b+1}{2}\sum
_{l=1}^{p-k_{w}-1} \gamma^l+\frac{b-1}{2}\gamma^{p-k_w},
\end{equation}
where
\[
k_w=\# \{ 1\le j\le p\dvtx |W_{w_j}(w|j-1)|=1 \}.
\]
This will imply that $\|F_p\|_\infty\le1+b/(1-\gamma)$ almost surely.

For the case $p=1$, that $ |\sum_{i=0}^{j}W_{i} |\le1+(b-1)\gamma/2$
for all $0\le j\le b-1$ is a direct consequence of our assumptions
$\mathbb{P} (\sum_{i=0}^{b-1}W_i=1 )=1$, and almost surely $\mbox
{either } |W_i|=1 \mbox{ and } (\sum_{k=0}^{i-1} W_i,\sum_{k=0}^i W_i )
\in\{(0,1),(1,0)\}, \mbox{ or } |W_i|\le\gamma$. Thus $\|F_1\|_\infty
\le1+(b-1)\gamma/2$.

Now suppose that $p\ge1$ and \eqref{fp} holds. Let $w\in\mathscr
{A}^p$. For every $1\le j\le b-1$, let $s_j(w)=\sum
_{i=0}^{j-1}W_{i}(w)$. By construction, we have
\begin{eqnarray*}
F_{p+1}(t_{wj})&=&F_p(t_w)+ \bigl(F_p(t_w+b^{-p})-F_p(t_w) \bigr)s_j(w).
\end{eqnarray*}
If $|W_j(w)|=1$, then by our assumption we have $k_{wj}=k_w+1$ and
$s_j(w)\in\{0,1\}$, so $F_{p+1}(t_{wj})\in\{F_p(t_w),F_p(t_w+b^{-p})\}
$ and \eqref{fp} holds. Otherwise, $k_{wj}=k_w$, and since $|Q(w)|\le
\gamma^{p-k_w}$ we get
\begin{eqnarray*}
|F_{p+1}(t_{wj})|&\le& |F_p(t_w)| +\gamma^{p-k_w} |s_j(w)|
\\
&\le&1+\frac{b+1}{2}\sum_{l=1}^{p-k_{w}-1} \gamma^l + \frac
{b-1}{2}\gamma^{p-k_w}+\gamma^{p-k_w}+ \frac{b-1}{2}\gamma^{p+1-k_w}
\\
&=&1+\frac{b+1}{2}\sum_{l=1}^{p-k_{wj}} \gamma^l + \frac{b-1}{2}\gamma
^{p+1-k_{wj}}.
\end{eqnarray*}

  Now we prove that the limit $F$ of $(F_n)_{n\ge1}$ is not
uniformly H\"older continuous. To see this, we can consider of any
$q>1$ the positive measure $\mu_q$ on $[0,1]$ obtained almost surely as
the derivative in the distribution sense of $F_{W^{(q)}}$ (see
Definition~\ref{Wbeta}). It follows from the equality $Q(w)=\Delta
F(I_w)$ and Theorem IV(i) of \cite{B1} that $\lim_{n\to\infty}\log(
|\Delta F(I_n(t))|)/-n\log(b)=\varphi_W'(q)$ for $\mu_q$-almost every~$t$
where $I_n(t)$ is the semi-open to the right $b$-adic interval of
generation $n$ containing $t$. Moreover, by construction we have $\lim
_{q\to\infty}\varphi_W'(q)=0$. This implies that, with probability 1,
there exists a sequence $(t_k)_{k\ge1}$ of points in $[0,1]$ such that
$\lim_{k\to\infty}h_F(t_k)=0$; hence $F$ is not uniformly H\"older.

\subsection[Proof of Theorem 2.3]{Proof of Theorem \protect\ref{timechange}}
The facts that $\beta\in(1,2)$ when $\mathbb{P}(\sum
_{i=0}^{b-1}W_i=1)<1$, as well as the properties of $W^{(\beta)}$ are immediate.

We will prove in the end of this section the following properties which
hold under the assumptions of Theorem~\ref{timechange}.
\begin{lemma}\label{osci}
Let $X$ and $X_{\beta}$ stand respectively for the oscillations of
$F_W$ and $F_{W^{(\beta)}}$ over $[0,1]$. For all $q\in\mathbb{R}$ we
have $\mathbb{E}(X^q)<\infty$ and $\mathbb{E}(X_\beta^q)<\infty$.
\end{lemma}

For every $w\in\mathscr{A}^{*}$, let $X(w)$ and $X_{\beta}(w)$ stand
respectively for the oscillations of $F^{[w]}_W$ and $F^{[w]}_{W^{(\beta
)}}$ over $[0,1]$. We deduce from the Lemma~\ref{osci} that, with
probability~1, for every $\varepsilon>0$, there exists $n_\varepsilon$
such that
%
\begin{equation}\label{osci1}
\forall n\ge n_\varepsilon,  \forall w\in\mathscr{A}^{n},  \forall
Y\in\{X,X_\beta\},\qquad b^{-n\varepsilon}\le Y(w)\le b^{n\varepsilon}.
\end{equation}
This implies that
\begin{eqnarray*}
b^{-n\varepsilon} |Q(w)|&\le& \operatorname{Osc}_{F_W}(I_w)=|Q(w)|X(w)\le
b^{n\varepsilon} |Q(w)|,
\\
b^{-n\varepsilon} |Q(w)|^\beta&\le& \operatorname{Osc}_{F_{W^{(\beta
)}}}(I_w)=|Q(w)|^\beta X_\beta(w)\le b^{n\varepsilon} |Q(w)|^\beta.
\end{eqnarray*}
Consequently,
\[
b^{-n(1+1/\beta) \varepsilon} \le\frac{\operatorname{Osc}_{F_W}(I_w)}{{\rm
Osc}_{F_{W^{(\beta)}}}(I_w)^{1/\beta}}=\frac{X(w)}{X_\beta(w)^{1/\beta
}}\le b^{n(1+1/\beta) \varepsilon}.
\]
Let $B=F_W\circ F_{W^{(\beta)}}^{-1}$. Let $J_w=F_{W^{(\beta)}}(I_w)$.
We have $|J_w|=\operatorname{Osc}_{F_{W^{(\beta)}}}(I_w)$, and ${\rm
Osc}_B(J_w)=\operatorname{Osc}_{F_W}(I_w)$, so the previous inequality is
equivalent to
\[
b^{-n(1+1/\beta) \varepsilon}\le \frac{\operatorname{Osc}_B(J_w)}{|J_w|^{1/\beta
}}=\frac{X(w)}{X_\beta(w)^{1/\beta}}\le b^{n(1+1/\beta) \varepsilon}.
\]
Under our assumptions, it is also true that [see Theorem
\ref{conv-1}(1)]
\[
\lim_{n\to\infty}\inf_{w\in\mathscr{A}^{n}} \frac{\log_b |J_w|}{-n}\ge
\alpha_0,
\]
where $\alpha_0=\sup_{p>0}\varphi_{W^{(\beta)}}(p)/p>0$ (in fact the
equality holds). Also, we have the following property (we postpone its
proof until after that of Lemma~\ref{osci}).

\begin{lemma}\label{ratio}
With probability 1, for every $\varepsilon>0$, there exists
$n_\varepsilon$ such that
%
\begin{equation}\label{ratio1}
\forall n\ge n_\varepsilon,\qquad  b^{-n\varepsilon} \le\inf_{w\in\mathscr
{A}^{n}} \inf_{0\le i\le b-1}\frac{|J_{wi}|}{|J_{w}|}\le1.
\end{equation}
\end{lemma}

We can also choose the random integer $n_\varepsilon$ so that for all
$n\ge n_\varepsilon$ \eqref{ratio1} holds as well as the property $
|J_w|\le b^{-n\alpha_0/2}$ for all $w\in\mathscr{A}^{n}$.

Let $t\in(0,1)$ and $0<r<\min_{w\in\mathscr{A}^{n_\varepsilon
+1}}|J_w|$. Let $w_2\in\mathscr{A}^{*}$ such that $|w_2|>n_\varepsilon
$, $t\in J_{w_2}\subset[t-r,t+r]$ and $|J_{w_2}|$ is maximal. Then let
$(w_1,w_3)\in\mathscr{A}^{*}\times\mathscr{A}^{*} $ such that $\min
(|w_1|,|w_3|)>n_\varepsilon$, $\min(|J_{w_1}|, |J_{w_3}|)\ge r$, the
intervals $J_{w_i}$, $i\in\{1,2,3\}$, are adjacent, $[t-r,t+r]\subset
J_{w_1}\cup J_{w_2}\cup J_{w_3}$, and $|J_{w_1}|+|J_{w_2}|+|J_{w_3}|$
is minimal. This constraint imposes that $|J_w|\le rb^{\varepsilon
|w|}$ for $w\in\{w_1,w_3\}$. Otherwise, due to \eqref{ratio1} we can
replace $J_w$ by one of its sons in the covering of $[t-r,t+r]$. Also,
we have $2r\ge|J_{w_2}|\ge r b^{-\varepsilon|w_2|}$. Otherwise, since
$t\in J_{w_2}$, due to \eqref{ratio1} we can replace $J_{w_2}$ by its
father, hence $|J_{w_2}|$ is not maximal.

Since
\[
\operatorname{Osc}_B(J_{w_2})\le\operatorname{Osc}_B([t-r,t+r])\le3\max_i \operatorname{Osc}_B(J_{w_i}),
\]
we have
\[
|J_{w_2}|^{1/\beta} b^{-|w_2|(1+1/\beta) \varepsilon}\le{\rm
Osc}_B([t-r,t+r])\le3 \max_i |J_{w_i}|^{1/\beta} b^{|w_i|(1+1/\beta)
\varepsilon}.
\]
Now we specify $\varepsilon<\alpha_0/4$. We can deduce from the
constraints on the length of the intervals $J_{w_i}$ that
$b^{\varepsilon|w_i|}\le r^{-4\varepsilon/\alpha_0}$. Consequently,
there exists a constant $C$ depending on $W$ only such that for $r$
small enough,
\[
r^{1/\beta+C\varepsilon}\le\operatorname{Osc}_B([t-r,t+r])\le r^{1/\beta
-C\varepsilon}.
\]
Since this holds for all $0<\varepsilon<\alpha_0/4$, almost surely for
all $t\in(0,1)$, we have in fact that, with probability 1, for all
$t\in(0,1)$, $\lim_{r\to0^+} \frac{\log( \operatorname{Osc}_B([t-r,t+r]))}{\log
(r)}=1/\beta$, hence $h_B(t)=1/\beta$.

\begin{pf*}{Proof of Lemma \protect\ref{osci}} The part concerning the
moments of positive orders is a consequence of Theorem \ref{conv-1}(1)
and the inequality $\operatorname{Osc}_f([0,1])\le2\|f\|_\infty$, which holds
for any continuous function $f$ on $[0,1]$.

For the moments of negative orders, the case of $X_\beta$, which is the
increment between 0 and 1 of the increasing function $F_{W^{(\beta)}}$,
is treated for instance in any of \cite{Mol,B1,Liu2001}. For $X$, we
just remark that we have
%
\begin{equation}\label{momentneg}
X\ge b^{-1}\sum_{k=0}^{b-1}|W_k|X(k).
\end{equation}
Moreover, the event $\{X=0\}$ is measurable with respect to $\bigcap
_{n\ge1} \sigma(W(w)\dvtx |w|\ge n)$ because the components of $W$ do not
vanish. Thus this event has probability 0 or 1. Since the function
$F_W$ is not almost surely equal to 0, $X$ is positive with probability
1. We can then use inequality \eqref{momentneg} in the same way as in
\cite{HoWa,Mol,B1} when $W$ is positive to prove that all the moments
of negative order of $X$ are finite as soon as the same property holds
for the random variables $|W_i|$.
\end{pf*}

\begin{pf*}{Proof of Lemma \protect\ref{ratio}} Let $(\varepsilon_k)_{k\ge
1}$ be a positive sequence converging to 0 at $\infty$. Since
$|J_w|=|Q(w)|^\beta X_\beta(w)$, due to \eqref{osci1}, for any $k\ge
1$, with probability~1, for $n$ large enough we have $b^{-n\varepsilon
_k} \le\inf_{w\in\mathscr{A}^{n}} \inf_{0\le i\le b-1}|W_{i}(w)|$.
This is an immediate consequence of the fact that all the moments of
negative order of the random variables $|W_i|$ are finite. Since the
set $\{\varepsilon_k\dvtx k\ge1\}$ is countable, we have the conclusion.
\end{pf*}

\subsection[Proof of Theorem 2.4]{Proof of Theorem~\protect\ref{notimechange}} Suppose that there
exist a continuous and increasing function $G$ defined on $[0,1]$ with
$G(0)=0$ as well as a monofractal continuous function $B$ defined on
$[G(0),G(1)]$ such that $F=B\circ G$. We denote by $H$ the H\"older
exponent of $B$ (notice that it may be random).

At first, suppose that $H\in(0,1]$. For every $\alpha\in(0,H)$ the
function $B$ is uniformly $\alpha$-H\"older, so there exists $C>0$ such
that $\operatorname{Osc}_B (I)\le C|I|^\alpha$. Consequently, for every $n\in
\mathbb{N}_+$ and $q\ge0$, we have
\[
\sum_{w\in\mathscr{A}^n}\operatorname{Osc}_F(I_w)^q\le
C^q\sum_{w\in\mathscr
{A}^n} |G(I_w)|^{\alpha q}.
\]
Since $G$ is increasing, taking $q=1/\alpha$ yields $
\sum_{w\in\mathcal{A}^n}\operatorname{Osc}_F(I_w)^{1/\alpha}\le C^{1/\alpha}
G(1)$, hence $\tau_F(1/\alpha)\ge0$. This is in contradiction to the
fact (established in \cite{BJ}) that $\tau_F=\varphi_W$ over $\mathbb
{R}_+$ almost surely.

Now we suppose that $H=0$. If $\dim ([0,1]\setminus\overline
E_G(\infty) )>0$, then we have $h_{F_W}(t)=0=h_B(G(t))$ at each $t\in
[0,1]\setminus\overline E_G(\infty)$. But since $\tau_F\ge\varphi_W$
over $\mathbb{R}_+$, we have $\dim E_F(0)\le\tau_F^*(0)\le\inf_{q\ge
0}-\varphi_W(q)=0$. This yields a contradiction.

\section[Proofs of Theorems 2.5, 2.6, 2.7 and 2.8]{Proofs of Theorems \protect\ref{div},
\protect\ref{div-caspart}, \protect\ref{div2}
and \protect\ref{div3}}\label{secdiv}

\subsection[Proof of Theorem 2.5]{Proof of Theorem \protect\ref{div}}
(1) If $\alpha\in\mathbb{R}$ and $n\ge1$ we define $F_{\alpha
,n}=b^{n\alpha}F_n$.

Let $\alpha_0=\varphi_W(p_0)/p_0=\varphi_W'(p_0)$, and
$W^{(p_0)}=b^{\varphi_W(p_0)}(|W_0|^{p_0},\dots,|W_{k-1}|^{p_0})$. By
construction, we have (see Definition~\ref{Wbeta}) $\varphi
_{W^{(p_0)}}(1)=\varphi_{W^{(p_0)}}'(1)=0$. We know from \cite{KP,DL}
that this implies that $\lim_{n\to\infty}\sum_{w\in\mathscr{A}^{n}}
\Delta F_{W^{(p_0)},n}(I_w)=0$ almost surely. Consequently, since
\[
\|F_{\alpha_0,n}\|_\infty\le \sum_{w\in\mathscr{A}^{n}} |\Delta
F_{\alpha_0,n}(I_w)|=\sum_{w\in\mathscr{A}^{n}} \Delta
F_{W^{(p_0)},n}(I_w)^{1/p_0}
\]
and $1/p_0\ge1$, we can conclude that $F_{\alpha,n}$ converges almost
surely uniformly to 0 for all $\alpha\le\alpha_0$.

Now let $\alpha>\alpha_0$ and set $V_\alpha=b^\alpha W$. We have
$\varphi_{V_\alpha}(p)=\varphi_W(p)-\alpha p$. Thus, $\varphi_{V_\alpha
}'(p_0)p_0-\varphi_{V_\alpha}(p_0)=\varphi_{W}'(p_0)p_0-\varphi
_{W}(p_0)=0$, $\varphi_{V_\alpha}(p_0)<0$ and \mbox{$\varphi_{V_\alpha
}'(p_0)<0$}. Consequently, we can find $p\in(0,p_0)$ such that $\varphi
_{V_\alpha}'(p)<0$ and $\varphi_{V_\alpha}'(p)p-\varphi_{V_\alpha
}(p)>0$. Let $V_\alpha^{(p)}=b^{\varphi_{V_\alpha}(p)}(|V_{\alpha
,0}|^p,\dots, V_{\alpha,b-1}|^p)$. We have $\varphi_{V_\alpha
^{(p)}}(1)=0$ and $\varphi_{V_\alpha}'(p)p-\varphi_{V_\alpha}(p)>0$ is
equivalent to $\varphi_{V_\alpha^{(p)}}'(1)>0$, so the nondecreasing
function~$F_{V_\alpha^{(p)}}$ is well defined and differs from 0 with
positive probability by Theorem~\ref{conv-1}. Let us denote by $\mu$
the measure $F_{V_\alpha^{(p)}}'$. It follows from Theorem IV(i) in~\cite{B1}
that, with probability~1, conditionally on $\mathscr{V}^c$,
we have $\mu\neq0$ and
\[
\lim_{n\to\infty}\frac{\log_b |\Delta F_{\alpha,n}(I_n(t))|}{-n}=-\sum
_{k=0}^{b-1}\mathbb{E}\bigl(V_{\alpha,k}^{(p)}\log_b|V_{\alpha,k}|\bigr)=\varphi
_{V_\alpha}'(p)
\qquad\mbox{for $\mu$-a.e. $t$},
\]
where $I_n(t)$ stands for the semi-open to the right $b$-adic interval
of generation $n$ containing $t$. Since $\varphi_{V_\alpha}'(p)<0$, we
conclude that $F_{\alpha,n}$ is unbounded.

(2) It is the same proof as in (1).

(3) We recall that in the conservative case we have $\mathbb
{P}(\mathscr{V}^c)=1$.
If $p_0<\infty$, the unboundedness result is proven as in 1.
If $p_0=\infty$ then $\sum_{i=0}^{b-1}\mathbb{P}(|W_i|=1)>1$ or $\mathbb
{P}(\#\{i\dvtx |W_i|=1\}=1)=1$.

We claim that, with probability 1, there exists an infinite word
$w_1\cdots w_n\cdots$ and an integer $n_0\ge1$ such that $|\Delta
F_{n_0}(I_{w_1\cdots w_{n_0}})|\neq0$ and $|W_{w_{n+1}}(w_1\cdots
w_n)|=1$ for all $n\ge n_0$. This implies that the absolute value of
the increment of $F_n$ over the interval $I_{w_1\cdots w_n}$ is equal
to the constant $|\Delta F_{n_0}(I_{w_1\cdots w_{n_0}})|$ for all $n\ge
n_0$, so that $(F_n)_{n\ge1}$ cannot converge uniformly.

Our claim is immediate if $\mathbb{P}(\#\{i\dvtx |W_i|=1\}=1)=1$. Now,
suppose that $\sum_{i=0}^{b-1}\mathbb{P}(|W_i|=1)>1$. If $n\ge1$ and
$w\in\mathscr{A}^n$, the set of words $v\in\bigcup_{p\ge n+1}\mathscr
{A}^p$ with prefix $w$ (i.e., $v|n=w$) such that $|W_{v_{k+1}}(v|k)|=1$
for all $n\le k\le|v|-1$ for a Galton--Watson tree $T(w)$ whose
offspring distribution is given by that of $N_1=\#\{i\dvtx |W_i|=1\}$.
Since $\sum_{i=0}^{b-1}\mathbb{P}(|W_i|=1)>1$ this tree is
supercritical, so it is finite with a probability $q_1<1$. Moreover,
the trees $T(w)$, $w\in\mathscr{A}^n$, are independent and independent
of $\mathcal{F}_n=\sigma(Q(w)\dvtx w\in\mathscr{A}^n)$.

Let $E_n$ be the event $\{\forall w\in\mathscr{A}^n,  Q(w)=0\mbox{
or } T(w)\mbox{ is finite}\}$ and $\mathscr{Z}_n=\{
w\in\mathscr{A}^n\dvtx
Q(w)\neq0\}$. We have $E_n=\{\forall w\in\mathscr{Z}_n,  T(w)\mbox{
is finite}\}$. Thus
\[
\mathbb{P}(E_n)= \mathbb{E}(\mathbb{P}(E_n|\mathcal{F}_n))=\mathbb
{E}(q_1^{\#\mathscr{Z}_n})=(\mathbb{E}(q_1^N))^n,
\]
where $N=\#\{0\le i\le b-1\dvtx |W_i|>0\}$ [see also Remark~\ref{Fn0}(3)].
Now, since $N\ge1$ almost surely and $q_1<1$, we have $\sum_{n\ge1}
\mathbb{P}(E_n)<\infty$. The conclusion follows from the
Borel--Cantelli Lemma.

The result about the nonexistence of a normalization making the
sequence weakly convergent is an obvious consequence of the fact that
we have $\Delta F_n(I_w)=\Delta F_{|w|}(I_w)$ for all $w\in\mathscr
{A}^{*}$ and $n\ge|w|$.

(4) follows from Theorem \ref{div-caspart}(1) and (2). Indeed,
suppose that $\alpha>\varphi_W(2)/2$ and $(b^{n\alpha} F_n)_{n\ge1}$
is bounded with positive probability. Since the boundedness of this
sequence is clearly an event which is measurable with respect to
$\bigcap_{p\ge1}\sigma(W(w), |w|\ge p)$, it occurs with probability
1. In this case, $X_n =r_n F_n$ tends almost surely uniformly to 0, as
$n\to\infty$. This contradicts the fact that the $L^2$ norm of $Z_n(1)$
converges to a positive value $\sigma$, and $X_n(1)$ is bounded in
$L^p$ norm for $p$ close to~$2^+$.

\subsection[Proof of Theorem 2.6]{Proof of Theorem \protect\ref{div-caspart}}
(1) Let $Y_0(w)=1$ and $Y_n(w)=F^{[w]}_n(1)$ for all $n\ge1$ and $w\in
\mathscr{A}^{*}$. Also, when $\varphi_W(2)<0$ let $\ell$ be the unique
solution of $\ell=b^{-\varphi_W(2)}\ell+\sum_{i\neq j}\mathbb
{E}(W_i\overline{W_j})$, that is, $\ell=\sum_{i\neq j}\mathbb
{E}(W_i\overline{W_j})/(1-b^{-\varphi_W(2)})$. We can deduce from (\ref
{foncmand}) that
\[
\mathbb{E}(|Y_n|^2)=b^{-\varphi_W(2)}\mathbb{E}(|Y_{n-1}|^2) + \sum
_{i\neq j}\mathbb{E}(W_i\overline{W_j}).
\]
Consequently, defining $v_n= \mathbb{E}(|Y_n|^2)$, for $n\ge0$ we have
$v_n=\ell+(v_0-\ell) b^{-n\varphi_W(2)} $ if $\varphi_W(2)< 0$ and
$v_n=v_0+n \mathbb{E}(W_i\overline{W_j})$ if $\varphi_W(2)=0$. If
$\varphi_W(2)<0$, we have $1=v_0\neq\ell$; otherwise we must have $1=
\sum_{i\neq j}\mathbb{E}(W_i\overline{W_j})+\break\mathbb{E} (\sum
_{i=0}^{b-1}|W_i|^2 )=\mathbb{E} (|\sum_{i=0}^{b-1}W_i|^2 )$, and since
we have $\mathbb{E} (|\sum_{i=0}^{b-1}W_i|^2 )^{1/2}\ge\break\mathbb{E} (
|\sum_{i=0}^{b-1}W_i| )\ge \mathbb{E} (\sum_{i=0}^{b-1}W_i )=1$, this
implies that we are in the conservative case. Also, if $\varphi
_W(2)=0$, then $\sum_{i\neq j}\mathbb{E}(W_i\overline{W_j})\neq0$;
otherwise $\mathbb{E} (|\sum_{i=0}^{b-1}W_i|^2 )= \mathbb{E} (\sum
_{i=0}^{b-1}|W_i|^2 )=b^{-\varphi_W(2)}=1$, and we have the same
contradiction as in the previous discussion.

Consequently, we have $\mathbb{E}(|Y_n|^2)\sim(1-\ell)b^{-n\varphi
_W(2)} =\sigma^2b^{-n\varphi_W(2)}$ if $\varphi_W(2)<0$ and $\mathbb
{E}(|Y_n|^2)\sim n\sum_{i\neq j}\mathbb{E}(W_i\overline{W_j})=\sigma^2
n$ if $\varphi_W(2)=0$.

(2) We denote by $\sigma_{n}$ the equivalent of $\sqrt{\mathbb
{E}(|F_n(1)|^2)}$ obtained in (1), that is, $\sigma_{n}=\sigma
b^{-n\varphi_W(2)/2}$ if $\varphi_W(2)<0$ and $\sigma_{n}=\sigma\sqrt
{n}$ if $\varphi_W(2)=0$, and we consider $Z_n=F_n/\mathscr{A}^{n}$
rather than $F_n/\sqrt{\mathbb{E}(|F_n(1)|^2)}$. For $w\in\mathscr
{A}^{*}$, we also denote $F_n^{[w]}/\sigma_{n}$ by $Z_n^{[w]}$.

We leave the reader to check the following simple properties for
$m,n\ge1$ and $w\in\mathscr{A}^m$: If $n>m$, then
%
\begin{equation}\label{self-sim2}
\Delta Z_{n}(I_w)=
Q(w)\cdot
\cases{
b^{m\varphi_W(2)/2} Z_{n-m}^{[w]}(1),&\quad\mbox{if $\varphi_W(2)<0$},
\vspace*{2pt}\cr
 \sqrt{\dfrac{n-m}{n}}Z_{n-m}^{[w]}(1),&\quad\mbox{if $\varphi_W(2)=0$},
}
\end{equation}
and if $1\le n\le m$ then
%
\begin{equation}\label{self-sim3}
\Delta Z_{n}(I_w)=Q(w|n)b^{n-m}\cdot
\cases{
b^{n\varphi_W(2)/2}/\sigma,&\quad\mbox{if $\varphi_W(2)<0$},
\vspace*{2pt}\cr
 1/\sigma\sqrt{n},&\quad\mbox{if $\varphi_W(2)=0$}.
}
\end{equation}

To simplify the notations, we denote $Z_n^{[w]}(1)$ by $\widetilde
Z_n(w)$. Also, we define $\widetilde W=b^{\varphi_W(2)/2}W$. By
construction, we have
%
\begin{equation}\label{self-simnorm}
\widetilde Z_{n+1}=
\cases{
\displaystyle\sum_{k=0}^{b-1}\widetilde W_k \widetilde Z_{n}(k),&\quad\mbox{if $\varphi
_W(2)<0$},
\vspace*{2pt}\cr
\displaystyle\sum_{k=0}^{b-1}\widetilde W_k \sqrt{\dfrac{n}{n+1}}\widetilde
Z_{n}(k),&\quad\mbox{if $\varphi_W(2)=0$}.
}
\end{equation}
This is formally the same [or almost the same when $\varphi_W(2)=0$]
equality as $F_{n+1}(1)=\sum_{k=0}^{b-1}W_k F_{n}^{[k]}(1)$, with the
same properties of independence and equidistribution. Moreover, we have
$\varphi_{\widetilde W}(p)=\varphi_W(p)-p\varphi_W(2)/2$, and
$\widetilde Z_n=\widetilde Z_n(\varnothing)$ is bounded in $L^2$ norm.
Consequently, when $\varphi_W(p)-p\varphi_W(2)/2>0$, the boundedness of
$\widetilde Z_n$ in $L^p$ is obtained by induction on the integer part
of $p$ like for the proof of the boundedness in $L^p$ of
$(F_n(1))_{n\ge1}$ when $\varphi_W(p)>0$ (see, for instance, \cite{KP,DL,B3}).

We now study the tightness of the sequence $(Z_n)_{n\ge1}$.
By Theorem~7.3 of \cite{Billingsley}, since $Z_n(0)=0$ almost surely for
all $n\ge1$, it is enough to show that for each positive $\varepsilon$
%
\begin{equation} \label{tightness}
\lim_{\delta\to0}\limsup_{n\to\infty} \mathbb{P}\bigl(\omega(
{Z}_n,\delta)> \varepsilon\bigr)=0
\end{equation}
(recall that $\omega(f,\cdot)$ stands for the modulus of continuity of
$f$ when $f\in\mathcal{C}([0,1])$).

We fix $p>2$ such that $\varphi_W(p)/p-\varphi_W(2)/2>0$ as well as
$\gamma>0$, and we estimate the $L^p$ norm of the $b$-adic increments
of $Z_n$ by using \eqref{self-sim2} and \eqref{self-sim3}.

If $n>m$, then
\begin{eqnarray*}
\sum_{w\in\mathscr{A}^m}\|\Delta Z_{n}(I_w)\|_p^p\le b^{-p m(\varphi
_W(p)/p-\varphi_W(2)/2)}\sup_{k\ge1}\|Z_{k}(1)\|_p^p,
\end{eqnarray*}
and if $1\le n\le m$, then
\begin{eqnarray*}
\sum_{w\in\mathscr{A}^m}\|\Delta Z_{n}(I_w)\|_p^p&\le& \sum_{w\in
\mathscr{A}^m} \mathbb{E}(|Q(w|n)|^p)b^{(n-m)p}
b^{np\varphi_W(2)/2}/\sigma^p
\\
&=& \sum_{v\in\mathscr{A}^{n}} b^{m-n}\mathbb{E}(|Q(v)|^p)b^{(n-m)p}
b^{np\varphi_W(2)/2)}/\sigma^p
\\
&=&b^{(m-n)(1-p)} b^{-np (\varphi_W(p)/p-\varphi_W(2)/2}/\sigma^p.
\end{eqnarray*}
Consequently, for all $m_0\ge1$ and $n\ge m_0$
\begin{eqnarray*}
&&\mathbb{P}(\exists m\ge m_0, \exists w\in\mathscr{A}^m\dvtx  | \Delta
Z_n(I_w)|>b^{-\gamma m} )
\\
&&\qquad \le \sum_{m\ge m_0} b^{p\gamma m} \sum_{w\in\mathscr{A}^m}\|\Delta
Z_{n}(I_w)\|_p^p
\\
&&\qquad \le \sup_{k\ge1}\|Z_{k}(1)\|_p^p \sum_{m=m_0}^{n-1} b^{-p m(\varphi
_W(p)/p-\varphi_W(2)/2-\gamma)}
\\
&& \qquad\quad {}+ \sigma^{-p} b^{-np (\varphi_W(p)/p-\varphi_W(2)/2-\gamma
)}\sum_{m\ge n}b^{(m-n)(1-p+\gamma p)}
\\
&&\qquad = O\bigl(b^{-pm_0(\varphi_W(p)/p-\varphi_W(2)/2-\gamma)}\bigr),
\end{eqnarray*}
if we choose $\gamma\in (0,\min(\varphi_W(p)/p-\varphi_W(2)/2,
(p-1)/p) )$. We fix such a $\gamma$, define $\alpha=\varphi
_W(p)/p-\varphi_W(2)/2-\gamma$ and notice thanks to \eqref{holder'}
that on the complement of $\{\exists m\ge m_0, \exists w\in\mathscr
{A}^m\dvtx \Delta Z_n(I_w)|>b^{-\gamma m} )\}$ we have $\omega
(Z_n,b^{-m_0})\le2b\cdot b^{-\gamma m_0 }/(1-b^{-\gamma})$.
Consequently, we have obtained that, for $\delta\in(0,1)$, and $m_0$
such that $b^{-m_0-1}< \delta\le b^{-m_0}$, we have
\[
\sup_{n\ge-\log_b(\delta)}\mathbb{P} \biggl(\omega(F_n, \delta)> \frac
{2b^{1+\gamma}}{1-b^{-\gamma}} \delta^\gamma\biggr)=O(\delta^\alpha).
\]

(3) The properties of $W^{(2)}$ are obvious. Suppose that $(Z_n)_{n\ge
1}$ converges in distribution to a continuous process $Z$, as $n\to
\infty$. Then the same holds for all the sequences $(Z^{[w]}_n)_{n\ge
1}$, and by using \eqref{self-simnorm} we see that the limit in
distribution of the $Z(1)$ must satisfy
\[
Z(1)=\sum_{k=0}^{b-1}\widetilde W_k Z^{[k]}(1).
\]
Moreover, $\mathbb{E}(Z(1))=0$ and $\mathbb{E}(Z^2(1))=1$.
Consequently, it follows from Theorem 4(ii) in \cite{Rosler} that the
characteristic function of $Z(1)$ is given by\vspace*{1pt} $\mathbb{E} (\exp
(-t^2\times F_{W^{(2)}}(1)/2) )$. It is clear that this is also the
characteristic function of $\widetilde Z= B\circ F_{W^{(2)}}(1)$.

Now let us prove that for all $p\ge1$, the vector $V_{p,n}=(\Delta
Z_{n}(I_w))_{w\in\mathscr{A}^p}$ converges in distribution to $(\Delta
B\circ F_{W^{(2)}}(I_w))_{w\in\mathscr{A}^p}$, as $n\to\infty$. This
will yield the conclusion.

Fix $p\ge1$ an integer. By definition of the processes $Z^{[w]}_n$,
for $n>p$, we have
\[
V_{p,n}= \bigl(r_{p,n}\widetilde Q(w)\cdot Z^{[w]}_{n-p}(1)\bigr)_{w\in\mathscr{A}^p}
\]
with $r_{p,n}=1$ if $\varphi_W(2)<0$ and $r_{p,n}=\sqrt{(n-p)/n}$
otherwise, and $\widetilde Q(w)=\prod_{k=1}^p\widetilde
W_{w_k}(w|k-1)$. Let $\mathcal{F}_p=\sigma(Q(w), w\in\mathscr{A}^p)$.
The random variables $Z^{[w]}_{n-p}(1)$ are independent and independent
of $\mathcal{F}_p$. Moreover, they converge in distribution to
$\widetilde Z$. Thus, if $\xi=(\xi_w)_{w\in\mathscr{A}^p}\in\mathbb
{R}^{\mathscr{A}^p}$, we have
\[
\lim_{n\to\infty}\mathbb{E} (\exp(i\langle\xi|V_{p,n}\rangle) |
\mathcal{F}_p )=\prod_{w\in\mathscr{A}^p}\phi_{\widetilde Z}(\xi_w
\widetilde Q(w)),
\]
where $\phi_{\widetilde Z}$ is the characteristic function of
$\widetilde Z$. Consequently, there exists a family $\{
F^{[w]}_{W^{(2)}}(1)\}_{w\in\mathscr{A}^p}$ of independent copies of
$F_{W^{(2)}}(1)$, this family being also independent of $\mathcal
{F}_p$, such that
\[
\lim_{n\to\infty}\mathbb{E} (\exp(i\langle\xi|V_{p,n}\rangle))=\mathbb
{E}\biggl(\prod_{w\in\mathscr{A}^p}\exp\bigl(-\xi_w^2\widetilde
Q(w)^2F^{[w]}_{W^{(2)}}(1)/2\bigr)\biggr).
\]

To see that the right-hand side is the characteristic function of
$U_p=(\Delta B\circ F_{W^{(2)}}(I_w))_{w\in\mathscr{A}^p}$, let us
first define $Q^{(2)}(w)=\prod_{k=1}^p W^{(2)}_{w_k}(w|k-1)$. Notice
that $Q^{(2)}(w)=\widetilde Q(w)^2$. By definition, $B\circ
F_{W^{(2)}}(I_w)=B(F_{W^{(2)}}(t_w+b^{-p}))-B(F_{W^{(2)}}(t_w))$.
Consequently, if $\xi\in\mathbb{R}^{\mathscr{A}^p}$,
\begin{eqnarray*}
\mathbb{E} \bigl(\exp(i\langle\xi|U_p\rangle )| \sigma({\mathcal{F}}_p\dvtx p\ge
1) \bigr)&=&\prod_{w\in\mathscr{A}^p}\exp\bigl(-\xi_w^2\Delta F_{W^{(2)}}(I_w)/2
\bigr)\\
&=&\prod_{w\in\mathscr{A}^p}\exp\bigl(-\xi_w^2Q^{(2)}(w)
F^{[w]}_{W^{(2)}}(1)/2 \bigr).
\end{eqnarray*}
Since $(\widetilde Q(w)^2)_{w\in\mathscr{A}^p}=(Q^{(2)}(w))_{w\in
\mathscr{A}^p}$, we have the result.

\subsection[Proof of Theorem 2.7]{Proof of
Theorem \protect\ref{div2}} The following proposition will
be crucial. We postpone its proof until the end of the section.

\begin{proposition}\label{detmom}
The probability distribution of $\widetilde Z=B\circ F_{W^{(2)}}(1)$ is
symmetric and determined by its moments. Specifically, for every
positive integer $q$ define $S_q=\{\beta=(\beta_0,\dots,
\beta_{b-1})\in\mathbb{N}^b\dvtx 0\le\beta_0,\dots,\beta_{b-1}<q,
\sum_{k=0}^{b-1}\beta_k=q\}$ and $M^{(q)}=\mathbb{E}(\widetilde Z^q)$.
Also let $\widetilde W=b^{\varphi_W(2)/2}W$. We have $M^{(2)}=1$ and
for every integer $p\ge2$
%
\begin{equation}\label{eqmom}
M^{(2p)}=\bigl(1-b^{-\varphi_{\widetilde W}(2p)}\bigr)^{-1} \mathop{\sum_{\beta
\in
S_{2p} }}_{ \beta_k \equiv0 [2]}\gamma_{ \beta} \mathbb{E} \Biggl(\prod
_{k=0}^{b-1} \widetilde W_k^{\beta_k} \Biggr)M^{(\beta_k)},
\end{equation}
where $ \gamma_{\beta_0,\dots, \beta_{b-1}} =\frac
{q!}{(\beta_0)!\cdots
(\beta_{b-1})!}$.
\end{proposition}

We now prove that the random variables $\widetilde Z_n=Z_n(1)$, $n\ge
1$, converge in distribution to $\widetilde Z$.

For all positive integers $q$ and $n$, let $\widetilde M_n^{(q)}=\mathbb
{E}(\widetilde Z_{n}^q)$, and as in Proposition \ref{detmom} let $S_q=\{
\beta=(\beta_0,\dots,
\beta_{b-1})\in\mathbb{N}^b\dvtx 0\le\beta_0,\dots,\beta_{b-1}<q,
\sum_{k=0}^{b-1}\beta_k=q\}$.

Due to Proposition \ref{detmom}, to get the convergence in distribution
of $\widetilde Z_n$ to $\widetilde Z$, it is enough to show the three
following properties:

\begin{longlist}[(1)]
\item[(1)] for every $p\ge0$ one has the property $(\mathcal{P}_{2p})$: $
\widetilde M^{(2p)}=\lim_{n\to\infty}\widetilde M_n^{(2p)}$ exists.
Moreover $\widetilde M^{(2)}=1$;

\item[(2)] for every $p\ge0$ one has the property $(\mathcal{P}_{2p+1})$: $
\lim_{n\to\infty}\widetilde M^{(2p+1)}_n=0$;

\item[(3)] for $p\ge2$
\[
\widetilde M^{(2p)}=\bigl(1-b^{-\varphi_{\widetilde W}(2p)}\bigr)^{-1} \mathop{\sum
_{\beta\in
S_{2p} }}_{ \beta_k \equiv0 [2]}\gamma_{ \beta} \mathbb{E} \Biggl(\prod
_{k=0}^{b-1} \widetilde W_k^{\beta_k} \Biggr)\widetilde M^{(\beta_k)}.
\]
\end{longlist}

For $n\ge1$ let $r_n=1$ if $\varphi_W(2)<0$ and $r_n=\sqrt{\frac
{n}{n+1}}$ otherwise. Let $q$ be an integer $\ge$3. Raising (\ref
{self-simnorm}) to
the power $q$ and taking the expectation yields
%
\begin{equation} \label{Mq}
\widetilde M_{n+1}^{(q)}=r_n^q b^{-\varphi_{\widetilde W}(q)}\widetilde
M_{n}^{(q)}+r_n^q\sum_{\beta\in S_q}\gamma_\beta
\mathbb{E} \Biggl(\prod_{k=0}^{b-1} \widetilde W_k^{\beta_k} \Biggr)\widetilde
M_{n}^{(\beta_k)}.
\end{equation}
We show by induction that $ ((\mathcal{P}_{2p-1}),(\mathcal{P}_{2p})
)$ holds for $p\ge1$, and we deduce the relation (3).

At first, $ ((\mathcal{P}_{1}),(\mathcal{P}_{2}) )$ holds by construction.
Suppose that $ ((\mathcal{P}_{2k-1}),(\mathcal{P}_{2k}) )$ holds for
$1\le k\le p-1$, with $p\ge2$. In particular, $M_n^{(\beta_k)}$ goes
to 0
as $n$ goes to $\infty$ if $\beta_k$ is an odd integer belonging to $[1,
2p-3]$.
Every element of the set $S_{2p-1}$ must contain an odd
component. Due to our induction assumption, this implies that in
relation (\ref{Mq}), the term
\[
r_n^{2p-1}b^{-(2p-1)/2}\sum_{\beta\in
S_{2p-1}}\gamma_{\beta} \mathbb{E} \Biggl(\prod_{k=0}^{b-1} \widetilde
W_k^{\beta_k} \Biggr)\widetilde M_{n}^{(\beta_k)}
\]
on the right-hand side goes to 0 at
$\infty$. This yields
\begin{eqnarray*}
\widetilde M^{(2p-1)}_{n+1}=r_n^{2p-1}b^{-\varphi_{\widetilde
W}(2p-1)}\widetilde M_{n}^{(2p-1)}+o(1)
\end{eqnarray*}
as $n\to\infty$. Since we have $r_n^{2p-1}b^{-\varphi_{\widetilde
W}(2p-1)}\le b^{-\varphi_{\widetilde W}(2p-1)}<1$, this yields $\lim
_{n\to\infty} \widetilde M_n^{(2p-1)}=0$, that is to say $(\mathcal{P}_{2p-1})$.

Now, the same argument as above shows that on the right-hand side of $
\widetilde M_{n+1}^{(2p)}$, we have
\begin{eqnarray*}
\lim_{n\to\infty} r_n^{2p}\sum_{\beta\in S_{2p}}\gamma_{ \beta} \mathbb
{E} \Biggl(\prod_{k=0}^{b-1} \widetilde W_k^{\beta_k} \Biggr)\widetilde
M_{n}^{(\beta_k)}=\mathop{\sum_{ \beta\in
S_{2p} }}_{ \beta_k \equiv0 [2]}\gamma_{ \beta}\mathbb{E} \Biggl(\prod
_{k=0}^{b-1} \widetilde W_k^{\beta_k} \Biggr)\widetilde M^{(\beta_k)}.
\end{eqnarray*}
Denote by $L$ the right-hand side of the above relation. By using (\ref
{Mq}) we deduce from the previous
lines that
\begin{eqnarray*} \label{rel}
\widetilde M_{n+1}^{(2p)}=r_n^{2p}b^{-\varphi_{\widetilde W}(2p)}
\widetilde M_{n}^{(2p)}+L+o(1).
\end{eqnarray*}
Consequently, since $\lim_{n\to\infty}r_n=1$, $\widetilde M_{n}^{(2p)}$
converges to the unique solution of $m=b^{-\varphi_{\widetilde W}(2p)}
m+L$. This yields both $(\mathcal{P}_{2p})$ and (3).

\begin{pf*}{Proof of Proposition \protect\ref{detmom}} Due to our assumption $\varphi
_{W^{(2)}}(p)>0$ for all $p>1$, it follows from Theorem 4.1 in \cite
{Liu1996} that $\limsup_{k\to\infty}\| F_{W^{(2)}}(1)\|_k/k<\infty$.
We also have $\limsup_{k\to\infty}\| B(1)\|_k/k<\infty$. Consequently,
since conditionally on $F_{W^{(2)}}$ we have $B\circ
F_{W^{(2)}}(1)\equiv F_{W^{(2)}} (1)^{1/2}B(1)$, we have $\limsup_{k\to
\infty}\| B\circ F_{W^{(2)}}(1)\|_k/k<\infty$. This ensures (see
Proposition 8.49 in \cite{Breiman}) that $B\circ F_{W^{(2)}}(1)$ is
determined by its moments.

Now we notice that there exists a vector $(\widetilde Z(0),\dots,
\widetilde Z(b-1))$ independent of $\widetilde W$ and whose components
are independent copies of $\widetilde Z$ such that
%
\begin{equation}\label{self-simfonc}
\widetilde Z\equiv\sum_{k=0}^{b-1}\widetilde W_k \widetilde Z(k).
\end{equation}
This is a consequence of Theorem 4(ii) in \cite{Rosler} as we said in
proving Theorem~\ref{div-caspart}(2). Since by construction $\mathbb
{E}(\widetilde Z)=0$ and $\mathbb{E}(\widetilde Z^2)=1$, we can exploit
(\ref{self-simfonc}) in the same way as (\ref{self-simnorm}) to get
(\ref{eqmom}).

For the sake of completeness, we give a direct argument for \eqref
{self-simfonc}. If we consider the random functions $F_{W^{(2)}}^{[0]},
\dots, F_{W^{(2)}}^{[0]}$ as well as $b$ independent Brownian motions
$B_0,\dots, B_{b-1}$, and if we set $\widetilde Z(k)=B_k\circ
F_{W^{(2)}}^{[k]}(1)$, then for $\xi\in\mathbb{R}$ we have
\begin{eqnarray*}
&&\mathbb{E} \Biggl(\exp\Biggl(i\Biggl\langle\xi\Big|\sum_{k=0}^{b-1}\widetilde W_k \widetilde
Z(k)\Biggr\rangle\Biggr) \Big| \sigma({\mathcal{F}}_p\dvtx p\ge1) \Biggr)
\\
&&\qquad =\prod_{k=0}^{b-1}\exp
\bigl(-\xi^2\widetilde W_k^2 F^{[k]}_{W^{(2)}}(1)/2 \bigr)
\\
&&\qquad =\prod
_{k=0}^{b-1}\exp\bigl(-\xi^2 W_k^{(2)} F^{[k]}_{W^{(2)}}(1)/2 \bigr)
\\
&&\qquad =\exp
\bigl(-\xi^2F_{W^{(2)}}(1)/2 \bigr).
\end{eqnarray*}
We have seen in the previous proof that the expectation of the
right-hand side is the characteristic function of $\widetilde Z$. This
yields (\ref{self-simfonc}). The previous computation shows that we can
obtain Theorem \ref{div-caspart}(4) without using \cite{Rosler}.
\end{pf*}

\subsection[Proof of Theorem 2.8]{Proof of Theorem~\protect\ref{div3}} At first we notice that under
our assumptions we have $\varphi_W(2)>0$. Then, since we can write
$(F_n-F)(1)=\sum_{w\in\mathscr{A}^{n}}Q(w)(1-F_W^{[w]}(1))$ and the
terms $(1-F_W^{[w]}(1))$ are centered, independent and independent of
the $Q(w)$, it is not difficult to see that $\mathbb
{E}((F_n-F)(1)^2)=\sigma^2 b^{-n\varphi_W(2)}$ where $\sigma=\sqrt
{\mathbb{E} ((1-F_W(1))^2 )}>0$. Also, $\varphi_{W^{(2)}}(p)>0$ in a
neighborhood of $1^+$, so $F_{W^{(2)}}$ is nondegenerate.

For $w\in\mathscr{A}^{*}$ let $R^{[w]}_n=(F^{[w]}_n-F^{[w]})/\sigma
b^{-n\varphi_W(2)/2}$. For $m,n\ge1$ and $w\in\mathscr{A}^m$: If
$n>m$, then
%
\begin{equation}\label{self-sim5}
\Delta R_{n}(I_w)=
|Q(w)|\cdot b^{m\varphi_W(2)/2} R_{n-m}^{[w]}(1),
\end{equation}
and if $1\le n\le m$, then
\begin{eqnarray*}
\Delta R_{n}(I_w)=b^{n-m}b^{n\varphi_W(2)/2}Q(w|n)/\sigma-b^{n\varphi
_W(2)/2}Q(w)F^{[w]}(1)/\sigma,
\end{eqnarray*}
and since $\varphi_W(2)>0$,
%
\begin{equation}\label{self-sim6}
\qquad\ |\Delta R_{n}(I_w)|\le b^{n-m}b^{n\varphi_W(2)/2}|Q(w|n)|/\sigma
+b^{m\varphi_W(2)/2}|Q(w)|\bigl|F^{[w]}(1)\bigr|/\sigma.
\end{equation}
We deduce from \eqref{self-sim5} that
\[
R_n(1)=\sum_{i=0}^{b-1} b^{\varphi_W(2)/2}W_iR^{[i]}_{n-1}(1).
\]
Then the same arguments as in the proof of Theorem~\ref{div}(2) show
that $R_n(1)$ is bounded in $L^p$ norm for $p$ in a neighborhood of
$2^+$. Also, \eqref{self-sim5} and \eqref{self-sim6} can be used to
prove the tightness of the sequence $(\mathcal{L}(R_n))_{n\ge1}$ like
\eqref{self-sim2} and \eqref{self-sim3} where used to prove the
tightness of $(\mathcal{L}(Z_n))_{n\ge1}$ in
the proof of Theorem~\ref{div}(2).

Now let us prove that for all $p\ge1$, the vector $V_{p,n}=(\Delta
R_{n}(I_w))_{w\in\mathscr{A}^p}$ converges in distribution to $(\Delta
B\circ F_{W^{(2)}}(I_w))_{w\in\mathscr{A}^p}$, as $n\to\infty$. This
will yield the conclusion. We adopt the same notations as in the proof
of Theorem~\ref{div}(2).

By definition of the processes $R^{[w]}_n$, for every $p\ge1$ and
$n>p$, we have
\[
V_{p,n}= \bigl(\widetilde Q(w)\cdot R^{[w]}_{n-p}(1)\bigr)_{w\in\mathscr{A}^p}.
\]
Consequently, due to the proof of Theorem~\ref{div}(2), it only remains
to prove that $R_p(1)$ converges in distribution to $\widetilde Z(1)$.
To do this, we remark that it follows from the argument used to prove
Proposition 4.1 and its Corollary 4.3 in \cite{OW} [which can be
directly applied to $R_p(1)$ when $W$ has nonnegative i.i.d.
components] that conditionally on $\mathcal{F}_p$,
$R_p(1)$ converges in distribution to a centered normal law of standard
deviation $\sqrt{F_{W^{(2)}}(1)}$. This implies that $R_p(1)$ converges
in law to $\widetilde Z(1)$ as $p\to\infty$ [see the expression of the
characteristic function of $\widetilde Z(1)$ in the proof of
Theorem~\ref{div}(3)].

%

\printaddresses

\end{document}